\documentclass[final,11pt,a4paper,times]{elsarticle}
\usepackage{amssymb}
\usepackage{amsthm}
\usepackage{latexsym}
\usepackage{amsmath}
\usepackage{color}
\usepackage{graphicx}
\usepackage{indentfirst}
\usepackage{mathrsfs}
\usepackage{url}
\usepackage[colorlinks,citecolor=blue,linkcolor=blue]{hyperref}
\usepackage{lineno,hyperref}
\allowdisplaybreaks
\newtheorem{theo}{Theorem}[section]

\newtheorem{lemm}[theo]{Lemma}
\newtheorem{prop}[theo]{Proposition}
\newtheorem{re}[theo]{Remark}

\newtheorem{definition}[theo]{Definition}

\numberwithin{equation}{section}

\journal{A}

\begin{document}
	  \begin{frontmatter}
		\title{ Quasi-periodic Swing via Weak KAM  Theory}
		\author {{Xun Niu $^{a}$ \footnote{E-mail address : niuxun21@mails.jlu.edu.cn; niuxunmath@163.com}} {Kaizhi Wang $^{b}$ \footnote{E-mail address : kzwang@sjtu.edu.cn}
Yong Li $^{a,c}$ \footnote{ E-mail address : liyong@jlu.edu.cn} }\\
{$^{a}$School of Mathematics and Statistics, Center for Mathematics and Interdisciplinary Sciences, Northeast Normal University,}{ Changchun, Jilin 130024,  China}\\
 {$^{b}$School of Mathematical Sciences, CMA-Shanghai, Shanghai Jiao Tong University, Shanghai 200240, China}\\
			    {$^{c}$Institute of Mathematics, Jilin University, Changchun 130012, China}\\
			    }
	    \cortext[mycorrespondingauthor]{Corresponding author}
		\begin{abstract}

          We investigate the dynamics of the quasi-periodic swing equations from the perspective of weak KAM theory.
          To this end, we firstly study a class of Hamiltonian systems.
          We obtain that the limit $u$, which derived from convergence of a sequence of functional minimizers,  satisfies  Hamilton-Jacobi equations
          in a weak sense.
          This is the so-called weak KAM solution.
          Meanwhile, we also get a  minimal measures $\mu$.
          Finally,  we derive that the existence of weak KAM solutions for the swing equations.
		\end{abstract}

		\begin{keyword}
			quasi-periodic swing \sep   weak KAM theory \sep minimal measures
		\end{keyword}

	\end{frontmatter}
	\tableofcontents
    \newpage
	\section{ Introduction \bf}\label{sec:1}
           Dynamical systems described  by Tonelli type Hamiltonians  have been extensively studied in recent years, and out of this came a most profound  result the weak KAM theory introduced by Fathi in \cite{F97,F98}.
           It builds a bridge between the complex dynamical behavior and the viscosity solutions theory of Hamilton-Jacobi equations, and provides another way to study invariant sets (such as Mather set, Aubry set and Ma\~{n}\'{e} sets).

          In terms of the class of Hamiltonians, just as far as we know,  Evans studied static Hamiltonian   via the weak KAM theory in \cite{E03,E04}.  Gomes  \cite{G} showed that a weak form of the KAM/Nekhoroshev theory for
         $H_\epsilon(p,x)=H_0(p)+\epsilon H_1(p,x)$.  For time-periodic Tonelli Lagrangian systems,  some interesting  results  was obtained in \cite{W-L,W-Y}  by using Lax-Oleinik semigroup. Recently, the papers \cite{N-J-L,N-W-L} studied the Hamiltonian  with normal variables and applied it to perturbation theory.

       In this paper, we mainly explore the coupled swing equation via the weak KAM theory.
          We establish the existence of weak KAM solutions for the corresponding Hamiltonian systems under periodic and quasi-periodic cases.

         The coupled swing equation plays a crucial role in the dynamics of power systems.
         The swing equation represents the motion of the rotor in the power systems, which controls the dynamics of  a synchronous generator.
           It is worth mentioning that the swing equation is also closely related to other systems such as mechanical pendulum mechanisms, Josephson junction circuits, phase-locked loops and so on, which can be referred to \cite{K}.
          Besides, for relevant work about this equation, please refer to \cite{A-S-V,G-S,Q-M-K-Z,S-M-V,S-G-H,S,S-M-H} and the references therein.

         Consider the following  coupled swing equation
         \begin{align}\label{Swing}
           \ddot{x}=\alpha+f(x,t),
         \end{align}
         where $x\in\mathbb{T}^n$, $\alpha:=(\alpha_1,\alpha_2,\cdots,\alpha_n)\in\mathbb{R}^n $ is  the mechanical power input, $\alpha_i\in [0,+\infty)$, $i=1,2,\cdots,n$,
 $f(x,t)=D_xV(x,t)$,
         \begin{align*}
          V(x,t):=-\sum_{i,j}^n\beta_{ij}(\omega t)\big(1-\cos(\lambda_i x_i+\lambda_jx_j)\big).
         \end{align*}
         Here, $\beta_{ij}(\cdot):\mathbb{R}^m\rightarrow\mathbb{R}$,
 $\omega:=(\omega_1,\omega_2,\cdots,\omega_m)\in\mathbb{R}^m$, $t\in\mathbb{R}$, $\lambda_i, \lambda_j $ are finite numbers.
         Thus, the corresponding Hamiltonian can be written as
         \begin{align}\label{Swing-H}
           H(x,y, t)=\frac{1}{2}\|y\|^2-\langle\alpha, x\rangle+\sum_{i,j}^n\beta_{ij}(\omega t
           )\big(1-\cos(\lambda_i x_i+\lambda_jx_j)\big),
         \end{align}
         where $y:=\dot{x}\in\mathbb{R}^n$.
         To study the dynamics of the above Hamiltonian system,
          we firstly provide the definition of the quasi-periodic function.
         \begin{definition}
         A continuous real function $f(t)$ is called quasi-periodic with frequency $\omega$, if there exists a rationally independent
real number sequence $\{\omega_l\}_{l=1}^m$ and a continuous function $F([\omega t]): \mathbb{T}^m\rightarrow\mathbb{R}$ such that
\begin{align*}
 f(t)=F([\omega t]),
\end{align*}
where $\omega:=(\omega_1,\omega_2,\cdots\omega_m)$, $[\cdot]$ denotes  the projection $\mathbb{R}$ onto $\mathbb{T}:=\mathbb{R}/\mathbb{Z}$.
         \end{definition}
         Based on  this definition, let $\phi_l:=[\omega_l t]$, $l=1,2,\cdots,m$,   Hamiltonian $H(x,y, t)$ then can be written in autonomous form,
          i.e. $H(x,y,\phi)$, which is thus $2\pi$ periodic in each components of $x$, $\phi$.
        Thus, the motion of Hamiltonian  $H(x,y,\phi)$ is described by the following equations
       \begin{equation}\label{1-3}
	      \left\{
				\begin{array}{l}
					\dot{x}=D_yH(x,y,\phi),\\
                    \dot{y}=-D_xH(x,y,\phi).
				\end{array}
			\right.
	  \end{equation}
      In partcularly, when $H$ is integrable, (\ref{1-3}) becomes
      \begin{equation*}\label{1-4}
	     \left\{
				\begin{array}{l}
					\dot{x}=D_yH(x,y,\phi),\\
                    \dot{y}=0.
				\end{array}
			\right.
	  \end{equation*}
      That is, it has an invariant torus $\mathbb{T}_{\omega}=\mathbb{T}^n\times\{y_0\}$ with frequency $\tilde{\omega}=D_yH(y_0)$, carrying a quasi-periodic flow $x(t)=\tilde{\omega} t+x_0$.
      However,  $H$ is non-integrable in general, such conclusion does not always hold true.

      Following the work of  Lions et al in \cite{L82,L83,L-P-V}, we can prove that, for each $P\in\mathbb{R}^n$, there exists a unique constant $\bar{H}(P)$ such that
        \begin{align}\label{1-5}
            H(x,P+D_xv(x,P,\phi),\phi)=\bar{H}(P)
         \end{align}
         has a viscosity solution $u(x,P,\phi)=x\cdot P+v(x,P,\phi)$.
         Then, we consider the viscosity solution $u(x,P,\phi)$ as a generating function, which implicitly define a coordinate transformation $(x,y)\rightarrow(X,P)$.
         Therefore, it can formally transform $H(x,y,\phi)$ into $\bar{H}(P)$, and (\ref{1-3}) becomes
         \begin{equation*}\label{1-A}
	       \left\{
				\begin{array}{l}
					\dot{X}=D_P\bar{H}(P),\\
                    \dot{P}=0.
				\end{array}
			\right.
	     \end{equation*}
          But the viscosity solution does not have sufficient smoothness, and hence it is  in general impossible to carry out this classical process.

      Inspired by Evans in \cite{E03}, we consider a variational problem to avoid this lack of regularity.
      We introduce the functional
      \begin{align*}
            I_k[v]:=\int_{\mathbb{T}^{n}\times\mathbb{T}^{m}} e^{kH(x,P+D_xv,\phi)}dxd\phi.
        \end{align*}
        Then, 
         we suppose throughout that the nonnegative Hamiltonian $H(x,y,\phi)$  is smooth with respect to $x$,~$y$, is continuous with respect to $\phi$, and satisfies the following hypotheses:
            \begin{itemize}
              \item (H1) Uniform convexity: There exists a constant $\gamma>0$ such that
                  \begin{align*}
                     \frac{\gamma}{2}\|\xi\|^2\leq H_{y_iy_j}\xi_i\xi_j,
                  \end{align*}
                  where $\|\cdot\|$ is  Euclidean norm, $\xi\in \mathbb{R}^n$, $i,j=1,2,\cdots n$;
              \item (H2) Growth bounded:
                  \begin{equation*}
                     \begin{array}{ccc}
                          |H_{x_ix_j}|\leq C(1+\|y\|^2), & |H_{x_iy_j}|\leq C(1+\|y\|), & |H_{y_iy_j}|\leq C,
                      \end{array}
                  \end{equation*}
                   where $|\cdot|$ stand for absolute value.
           \end{itemize}
           This implies
             \begin{align*}
                \begin{array}{ccc}
                   |H_{x_i}|\leq C(1+\|y\|^2), & |H_{y_i}|\leq C(1+\|y\|), & \frac{\gamma}{2}\|y\|^2-C\leq H\leq C(1+\|y\|^2).
                \end{array}
             \end{align*}

           According to the Proposition \ref{prop-5.1}, for each $k\in\mathbb{R}^+$,  we have a smooth solution $v^k$ (i.e. minimizer of $I_k[\cdot]$ ) for the Euler-Lagrange equation corresponding to the functional $I_k[\cdot]$.
       The sequence of minimizers $\{v^k\}_{k=1}^\infty$  indeed exists, and is convergent  by Lemma \ref{lem-3.1}.
         Therefore, a weak solution of the Hamilton-Jacobi equation (\ref{1-5}) is obtained,
         through a series of techniques,  so that the above process can be transformed in the sense of a minimal measure.
         This measure is closely related to the Mather set, for details see
         \cite{E02, E03,M,N-J-L}.
          Hence we call this  solution the weak KAM solution.
      \begin{theo}\label{th1}
	         Suppose that the  conditions $(H1)$, $(H2)$  hold.
          Then, there exists a solution $u(x,P,\phi):=x\cdot P+v(x,P,\phi)$ satisfies the Hamilton-Jacobi equation
            \begin{align*}
            H(x,P+D_xv,\phi)=\bar{H}(P)
         \end{align*}
             in the sense of  the minimal measure $\sigma$ and
             \begin{align*}
            {\rm div}_x[\sigma D_yH(x,D_xu,\phi)]=0
         \end{align*}
         in  $\mathbb{T}^{n}\times\mathbb{T}^{m}$.
          \end{theo}

           \begin{re}
            A weak KAM solution  is  also a viscosity solution of Aronsson's equation
             \begin{align*}
              H_{y_i}(x,D_xu,\phi)H_{y_j}(x,D_xu,\phi)u_{x_ix_j}+H_{x_i}(x,D_xu,\phi) H_{y_i}(x,D_xu,\phi)=0.
             \end{align*}
           \end{re}
         Finally, we come back to the dynamics of the swing equation.
         Through verification,  the Hamiltonian (\ref{Swing-H}) satisfies the conditions $(H1)$, $(H2)$.
          According to Theorem \ref{th1}, we obtain the existence of weak KAM solution of the Hamiltonian system corresponding to  the swing equation (\ref{Swing}), that is,
\begin{theo}\label{th2}
	       The Hamiltonian system corresponding to the swing equation
      \begin{align*}
        \ddot{x}=\alpha-D_x\bigg(-\sum_{i,j}^n\beta_{ij}(\omega t)\big(1-\cos(\lambda_i x_i+\lambda_jx_j)\big)\bigg)
      \end{align*}
      has a weak KAM solution.
          \end{theo}
      In addition, we consider the periodic coupled swing equation
      \begin{align}\label{Swing-P}
           \ddot{x}=f(x,t),
         \end{align}
         where $x\in\mathbb{T}^n$, $t\in\mathbb{R}^1$, $f(x,t)=D_xV(x,t)$,
         \begin{align*}
          V(x,t):=-\sum_{i,j}^n\beta_{ij}(t)\big(1-\cos( x_i+x_j)\big),
         \end{align*}
       and $\beta_{ij}(t)=\beta_{ij}(t+2\pi)$.
       Then,  similarly, we obtain the existence of weak KAM solutions in this case.
       Further, this weak KAM solution is a classical solution of an elliptic equation.
  \begin{theo}\label{th-P}
      The elliptic equation
             \begin{align*}
            \sigma u_{x_ix_i}+\sigma_{x_i}u_{x_i}=0
         \end{align*}
         has a periodic classical solution $u(x,P,\phi):=x\cdot P+v(x,P,\phi)$,
         where $\phi:=[t]\in\mathbb{T}^1$.
\end{theo}
     \textbf{Structure of the paper.}
      In Sect. \ref{sec:2}, we present the  outline of proof to the main results (Theorem \ref{th1} and Theorem \ref{th2}).
      In order to get Theorem \ref{th1}, we start with some estimates (Lemma \ref{lem-3.1}) and prove the tightness of the probability measure sequence (Proposition \ref{pro-3.2}) in Sect. \ref{sec:3}.
      Then, we obtain  the limit of function $u$ and  the limit  of measure $\mu$.
      Further, in Sect. \ref{sec:4}, we prove that the probability measure $\mu$ is a minimal measure (Proposition  \ref{pro-4.4}), and $u$ satisfies some Hamilton-Jacobi equation in the sense of  the minimal measure  (Theorem \ref{th3}).
      In Appendix, we discuss smooth solutions of a kind of divergence equation (Proposition \ref{prop-5.1}) by the continuation method, which is crucial to our approach.
       Besides, We  supplement the proofs of some lemmas.
          \section{ Proof of  main theorems \bf}\label{sec:2}
       Here,  we present the  proof framework of Theorem \ref{th1} and Theorem \ref{th2}.

       \subsection{\textbf{Proof of  Theorem \ref{th1}}}\label{2.2}
       The proof of Theorem \ref{th1} is quite complicated, and to give the reader a clearer sense of our idea, we only give the sketch of the proof in  this section.
       The detailed proof is explained in Sect.  \ref{sec:3} and Sect. \ref{sec:4}.
          \begin{proof}
          Motivated by the inf-max formula
        \begin{align*}
           \bar{H}(P)=\mathop {\inf}_{v\in C^1(\mathbb{T}^{n})\times C(\mathbb{T}^{m})}\mathop {\max}_{(x,\phi)\in\mathbb{T}^{n}\times\mathbb{T}^{m}}H(x,P+D_xv,\phi),
         \end{align*}
         we consider a variational problem and introduce a functional
        \begin{align}\label{1-6}
            I_k[v]=\int_{\mathbb{T}^{n}\times\mathbb{T}^{m}} e^{kH(x,P+D_xv,\phi)}dxd\phi
        \end{align}
        for $k\in\mathbb{N}^+$.
       Firstly, we  prove the existence of minimizers  for functional $I_k[\cdot]$ by the continuation method (Proposition \ref{prop-5.1}),
        and  call $v^k\in C^1(\mathbb{T}^{n})\times C(\mathbb{T}^{m})$ the minimizer of $I_k[\cdot]$.
         Further, the minimizer of $I_k[\cdot]$ is unique, once we claim that the condition
        \begin{align*}
           \int_{\mathbb{T}^{n}\times\mathbb{T}^{m}} v^k dxd\phi=0
        \end{align*}
          holds.
          Besides, the inf-max formula can refer to \cite{C-I-P-P,E03,N-J-L}.
         Thus, the corresponding Euler-Lagrange equation for (\ref{1-6}) is
        \begin{align}\label{1-7}
            {\rm div}_x[e^{k[H(x,P+D_xv^k,\phi)}D_yH(x,P+D_xv^k,\phi)]=0.
        \end{align}
        Let $u^k(x,P,\phi)=x\cdot P+v^k(x,P,\phi)$.
         We  construct the sequence of minimizers $\{u^k\}_{k=1}^{\infty}$,
        and obtain the limit of function $u(x,P,\phi)$ by estimating $\{u^k\}_{k=1}^{\infty}$ (Lemma \ref{lem-3.1}).

        Next, we define  probability measure
        \begin{align*}
             \sigma^k:=\frac{e^{kH(x,P+D_xv^k,\phi)}}{\int_{\mathbb{T}^{n}\times\mathbb{T}^{m}} e^{kH(x,P+D_xv^k,\phi)}dxd\phi}=e^{k[H(x,P+D_xv^k,\phi)-\bar{H}^k(P)]},
        \end{align*}
        where
        \begin{align*}
           \bar{H}^k(P)=\frac{1}{k}\log\bigg(\int_{\mathbb{T}^{n}\times\mathbb{T}^{m}} e^{kH(x,P+D_xv^k,\phi)}dxd\phi\bigg).
        \end{align*}
       Obviously, we can right away get that $\sigma^k\geq0$, $\int_{\mathbb{T}^{n}\times\mathbb{T}^{m}}d\sigma^k=1$, where  $d\sigma^k:=\sigma^kdxd\phi$.
       Then, Euler-Lagrange equation (\ref{1-7}) becomes
       \begin{align}\label{1-8}
           {\rm div}_x[\sigma^kD_yH(x,P+D_xv^k,\phi)]=0.
      \end{align}
       On the other hands, for lifting the measure space of $\sigma^k$  to $\mathbb{T}^n\times\mathbb{R}^n\times\mathbb{T}^m$, we define
       \begin{align*}
          \mu^k=: \delta_{\{\beta=D_yH(x,D_xu^k,\phi)\}}\sigma^k,
       \end{align*}
       that is, for each  bounded continuous function $\Psi$,
       \begin{align*}
         \int_{\mathbb{T}^n\times\mathbb{R}^n\times\mathbb{T}^m}\Psi(x,\beta,\phi)d\mu^k=\int_{\mathbb{T}^{n}\times\mathbb{T}^{m}}\Psi(x,D_yH(x,D_xu^k,\phi),\phi)d\sigma^k.
       \end{align*}
 We  construct the  the sequence  of probability measure $\{\mu^k\}_{k=1}^{\infty}$,
  and derive the  tightness of probability measure sequence $\{\mu^k\}_{k=1}^{\infty}$ (Proposition \ref{pro-3.2}).

  Finally, in Sect. \ref{sec:4}, we will show that $\mu$  is a  minimal measure (Proposition  \ref{pro-4.4}) and $u(x,P,\phi)$ satisfies the Hamilton-Jacobi  equation
       \begin{align}\label{3-C}
          H(x,D_xu,\phi)=\bar{H}(P)
       \end{align}
       $\sigma-$ almost everywhere (Theorem \ref{th3}), where $\sigma$ is  the projection of $\mu$ into $\mathbb{T}^{n}\times\mathbb{T}^{m}$.
       Therefore,  it can transform $H(x,y,\phi)$ into $\bar{H}(P)$ in the sense of  the minimal measure, and (\ref{1-3}) becomes
      \begin{equation*}
	  \left\{
				\begin{array}{l}
					\dot{X}=D_P\bar{H}(P),\\
                    \dot{P}=0.
				\end{array}
			\right.
	 \end{equation*}
     Therefore, the systems $H(x,y,\phi)$ has  weak KAM solutions $u(x,P,\phi)$.
     In addition, according to the continuous dependence of the solution on the parameters, the weak KAM solution $u(x,P,\phi)$ is continuous with respect to $\phi$ (Proposition \ref{pro-4.7}).
      \end{proof}
      \subsection{\textbf{Proof of  Theorem \ref{th2}}}\label{2.3}
      \begin{proof}
      According to the coupled swing equation (\ref{Swing}), we can obtain the corresponding explicit Hamiltonian as
      \begin{align*}
             H(x,y,\omega t)=\frac{1}{2}\|y\|^2-\langle\alpha, x
      \rangle+\sum_{i,j}^n\beta_{ij}([\omega t])\big(1-\cos(\lambda_i x_i+\lambda_jx_j)\big),
      \end{align*}
      where $x\in\mathbb{T}^n$, $y\in\mathbb{R}^n$, $t\in\mathbb{R}$,  $\omega:=(\omega_1,\omega_2,\cdots,\omega_m)\in\mathbb{R}^m$, $\alpha$ is an $n$-dimensional constant vector, $\lambda_i$, $\lambda_j$ are finite numbers.
Let $\phi_l:=\omega_l t$, $l=1,2,\cdots,m$, the above Hamiltonian $H(x,y,\omega t)$
can be written as
\begin{align*}
             H(x,y,\phi)=\frac{1}{2}\|y\|^2-\langle\alpha, x
      \rangle+\sum_{i,j}^n\beta_{ij}(\phi)\big(1-\cos(\lambda_i x_i+\lambda_jx_j)\big),
      \end{align*}
      where $\phi \in\mathbb{T}^m$.
      Then, we derive that
      \begin{equation*}
       \begin{array}{ccc}
       H_{y_iy_j}=\frac{1}{2}, &  H_{x_ix_j}=\sum_{i,j}^n\beta_{ij}(\phi)\lambda_i\lambda_j\cos(\lambda_i x_i+\lambda_jx_j), &
         H_{x_iy_j}=0.
       \end{array}
      \end{equation*}
      Hence the above  Hamiltonian $H(x,y,\phi)$  satisfies the following hypotheses:
            \begin{itemize}
              \item (H1) Uniform Convexity:
                  \begin{align*}
                     \frac{\gamma}{2}\|\xi\|^2\leq H_{y_iy_j}\xi_i\xi_j,
                  \end{align*}
                  where $\gamma\in(0,1]$, $\xi\in \mathbb{R}^n$, $i,j=1,2,\cdots n$;
              \item (H2) Growth Bounded:
                  \begin{equation*}
                     \begin{array}{ccc}
                          |H_{x_ix_j}|\leq C, & |H_{x_iy_j}|=0, & |H_{y_iy_j}|\leq \frac{1}{2},
                      \end{array}
                  \end{equation*}
           \end{itemize}
      Finally, according to Theorem \ref{th1}, we obtain that the swing equation
      \begin{align*}
        \ddot{x}=\alpha-D_x\bigg(-\sum_{i,j}^n\beta_{ij}([\omega t])\big(1-\cos(\lambda_i x_i+\lambda_jx_j)\big)\bigg)
      \end{align*}
      has a weak KAM solution.
      \end{proof}
      \subsection{\textbf{Proof of  Theorem \ref{th-P}}}\label{2.3}
      \begin{proof}
      This proof is divided into six steps.

      \noindent \textbf{Step 1.}
      According to the periodic coupled swing equation (\ref{th-P}),  the corresponding Hamiltonian can be written as
         \begin{align}\label{Swing-H-P}
           H(x,y, t)=\frac{1}{2}\|y\|^2+\sum_{i,j}^n\beta_{ij}(t)\big(1-\cos( x_i+x_j)\big),
         \end{align}
         where $y:=\dot{x}\in\mathbb{R}^n$.
         Let $\phi:=[t]\in\mathbb{T}^{1}$,
         then, Hamiltonian (\ref{Swing-H-P}) can be written
         \begin{align}\label{Swing-H-P2}
           H(x,y, t)=\frac{1}{2}\|y\|^2+\sum_{i,j}^n\beta_{ij}(\phi)\big(1-\cos( x_i+x_j)\big).
         \end{align}
      Therefore, based on Theorem \ref{th1}, there exists a solution $u(x,P,\phi):=x\cdot P+v(x,P,\phi)$ satisfies the Hamilton-Jacobi equation
            \begin{align*}
            H(x,P+D_xv,\phi)=\bar{H}(P)
         \end{align*}
             in the sense of  the minimal measure $\sigma$ and
             \begin{align*}
          \sigma u_{x_ix_i}+\sigma_{x_i} u_{x_i}=0
         \end{align*}
         in  $\mathbb{T}^{n}\times\mathbb{T}^{1}$.
        That is, the Hamiltonian (\ref{Swing-H-P2}) has a weak KAM solution.

        \noindent\textbf{Step 2.}
        Without loss of generality, we assume $\sigma^k\neq0$.
         Because $u^k$ is the smooth solution of
         \begin{align*}
           {\rm div}_x[e^{k[H(x,D_xu^k,\phi)}D_yH(x,D_xu^k,\phi)]=0,
         \end{align*}
         that is,
         \begin{align*}
           u_{x_ix_i}^k+\frac{1}{\sigma^k}\sigma^k_{x_i}u^k_{x_i}=0.
         \end{align*}
       Let $\mathcal{O}u^k:=u^k_{x_ix_i}+au^k$, where $a<0$.
        We observed that $\mathcal{O}u^k=0$ only has trivial solutions, that is,
        $u^k\equiv0$.
        According to the Fredholm alternative theorem,
        the operator $\mathcal{O}$ has a solution $u^k$
         for the equation $\mathcal{O}u^k=Nu^k$,
        and this solution can be expressed as
        \begin{align}\label{th3-3}
          u^k(x)=\int_{\mathbb{T}^n}G(x,z)Nu^k(z)dz,
        \end{align}
        where $G(x,z)$ is an integral kernel,
        $Nu^k:=-\frac{1}{\sigma^k}\sigma^k_{x_i}u^k_{x_i}+au^k$.
        According to Lemma \ref{lem-3.1} and Theorem \ref{th1}, there exists a solution $u(x,P,\phi):=x\cdot P+v(x,P,\phi)$ such that
        \begin{align*}
         u(x,P,\phi)=\lim_{k\to\infty}u^k(x)=\lim_{k\to\infty}\int_{\mathbb{T}^n}G(x,z)Nu^k(z)dz,
        \end{align*}
        and the solution $u(x,P,\phi)$  satisfies
        \begin{align*}
         \int_{\mathbb{T}^{n}\times\mathbb{T}^{m}} {\rm div}_x[\sigma D_yH(x,D_xu,\phi)] \varphi dxd\phi=0
        \end{align*}
        for all $\varphi\in C^1(\mathbb{T}^{n}\times\mathbb{T}^{m})$.
        That is, $u(x,P,\phi)$ is a weak solution of the equation ${\rm div}_x[\sigma D_yH(x,D_xu,\phi)]=0$.
         Because $u^k_{x_ix_i}+au^k=Nu^k$, we have

        \begin{align*}
          \begin{array}{c}
            u_{x_1x_1}^k+\frac{a}{n}u^k=\frac{1}{n}Nu^k, \\
         u_{x_2x_2}^k+\frac{a}{n}u^k=\frac{1}{n}Nu^k, \\
         \vdots\\
         u_{x_nx_n}^k+\frac{a}{n}u^k=\frac{1}{n}N_nu^k.
          \end{array}
        \end{align*}
        Then, based on Lagrange's method of variation of constant,
        we obtain  the Green's function corresponding to each equation, that is,
        \begin{align*}
          G_i(x_i, z_i)=\sqrt{\frac{n}{|a|}}\sin(\sqrt{\frac{|a|}{n}}|x_i-z_i|)
        \end{align*}
         where $i=1,2,\cdots,n$.
        Let $G(x,z):=\sum_{i=1}^nG_i(x_i,z_i)$,
        \begin{align}\label{th3-3}
          u^k(x)=\int_{\mathbb{T}^n}G(x,z)Nu^k(z)dz.
        \end{align}
Based on the expression of $u^k$ (\ref{th3-3}), we have
        \begin{align*}
          u^k_{x_i}(x)=\int_{\mathbb{T}^n}G(z)N_{x_i}u^k(x,z)dz,
        \end{align*}
        and
         \begin{align*}
          u^k_{x_ix_i}(x)=\int_{\mathbb{T}^n}G(z)N_{x_ix_i}u^k(x,z)dz,
        \end{align*}
        where
        \begin{align*}
          N_{x_i}u^k(x,z):=\bigg(\frac{1}{(\sigma^k)^2}
          (\sigma^k_{x_i})^2u^k_{x_i}-\frac{1}{\sigma^k}\sigma^k_{x_ix_i}u^k_{x_i}
          -\frac{1}{\sigma^k}u^k_{x_i}u^k_{x_ix_i}+au^k_{x_i}\bigg)(x,z),
        \end{align*}
        and
        \begin{align*}
          N_{x_ix_i}u^k(x,z)&:=\bigg(-\frac{2}{(\sigma^k)^3}
          (\sigma^k_{x_i})^3u^k_{x_i}+\frac{3}{(\sigma^k)^2}\sigma^k_{x_i}
          \sigma^k_{x_ix_i}u^k_{x_i}
          +\frac{1}{(\sigma^k)^2}(\sigma^k_{x_i})^2u^k_{x_ix_i}
          -\frac{1}{\sigma^k}\sigma^k_{x_ix_ix_i}u^k_{x_i}\\
          &-\frac{1}{\sigma^k}\sigma^k_{x_ix_i}u^k_{x_ix_i}
          +\frac{1}{(\sigma^k)^2}\sigma^k_{x_i}u^k_{x_i}u^k_{x_ix_i}
          -\frac{1}{\sigma^k}(u^k_{x_ix_i})^2
          -\frac{1}{\sigma^k}u^k_{x_i}u^k_{x_ix_ix_i}
          +au^k_{x_ix_i}\bigg)(x,z).
        \end{align*}

        \noindent\textbf{Step 3.}
        Next, we discuss the convergence of $\sigma^k_{x_i}$.
         We  observe the definition of
        \begin{align*}
             \sigma^k:=\frac{e^{kH(x,D_xu^k,\phi)}}{\int_{\mathbb{T}^{n}\times\mathbb{T}^{m}} e^{kH(x,D_xu^k,\phi)}dxd\phi}=e^{k[H(x,D_xu^k,\phi)-\bar{H}^k(P)]},
        \end{align*}
        and obtain
        \begin{align*}
          \sigma^k_{x_i}=k[H(x,D_xu^k,\phi)]_{x_i}\sigma^k.
        \end{align*}
        Let  $\sigma^{(1),k}:=\sigma^k_{x_i}$, $d\sigma^{(1),k}:=\sigma^{(1),k}dxd\phi$,
         $\mu^{(1),k}$ lift measure $\sigma^{(1),k}$ to space $\mathbb{T}^{n}\times\mathbb{R}^{n}\times\mathbb{T}^{1}$.
        Then,  for sufficiently large $M$,
         \begin{align*}
          &\mu^{(1),k}\{(x,\beta,\phi)~|~~\|\beta\|\geq M\}\\
          &\leq \sigma^{(1),k}\{(x,\phi)~|~H(x,D_xu^k,\phi)\geq \bar{H}^k(P)+1\}
          =\int_{H(x,D_xu^k,\phi)\geq \bar{H}^k(P)+1}d\sigma^{(1),k}\\
          &\leq\int_{\mathbb{T}^{n}\times\mathbb{T}^{1}}d\sigma^{(1),k}
          \leq\int_{\mathbb{T}^{n}\times\mathbb{T}^{1}}e^{\rho^{(1)}k
          (H(x,D_xu^k,\phi)-\bar{H}^k(P)-1)}d\sigma^{(1),k}\\
          &=e^{-\rho^{(1)}k}\int_{\mathbb{T}^{n}\times\mathbb{T}^{1}}
          (\sigma^k)^{\rho^{(1)}}k[H(x,D_xu^k,\phi)]_{x_i}\sigma^kdxd\phi\\
          &\leq e^{-\rho^{(1)}k}\bigg(\frac{k}{2}\int_{\mathbb{T}^{n}\times\mathbb{T}^{1}}
          (\sigma^k)^{2\rho^{(1)}+1}dxd\phi
          +\frac{1}{2}\int_{\mathbb{T}^{n}\times\mathbb{T}^{1}}k[H(x,D_xu^k,\phi)]^2_{x_i}d\sigma^k\bigg),
         \end{align*}
         where $\rho^{(1)}$  is a finite number and we will provide the specific form below.
         According to Proposition \ref{pro-3.2}, we have
         \begin{align*}
        \int_{\mathbb{T}^{n}\times\mathbb{T}^{1}} k [H(x,D_xu,\phi)]_{x_i}^2d\sigma
        \leq C,
        \end{align*}
        and
        \begin{align*}
        \int_{\mathbb{T}^{n}\times\mathbb{T}^{1}}\bigg(\|D_x(\sigma^{\frac{1}{2}})\|^2+(\sigma^{\frac{1}{2}})^2~~\bigg)dxd\phi
        \leq Ck.
        \end{align*}
         Based on Sobolev's embedded  inequality, that is, $||(\sigma^k)^{\frac{1}{2}}||_{L^{2(2\rho^{(1)}+1)}(\mathbb{T}^{n}\times\mathbb{T}^{1})} \leq C ||(\sigma^k)^{\frac{1}{2}}||_{W^{1,2}(\mathbb{T}^{n}\times\mathbb{T}^{1})}$, and
        \begin{align*}
            \rho^{(1)}=\left\{
                \begin{array}{lr}
                  \frac{1}{n-1}, & n\geq2; \\
                  {\rm any ~positive ~number}, & n=1.
                \end{array}
             \right.
        \end{align*}
         we have
        \begin{align*}
          \bigg(\int_{\mathbb{T}^{n}\times\mathbb{T}^{1}}
          (\sigma^k)^{2\rho^{(1)}+1}dxd\phi\bigg)^{\frac{1}{2\rho^{(1)}+1}}
          \leq\int_{\mathbb{T}^{n}\times\mathbb{T}^{1}}\bigg(\|D_x(\sigma^{\frac{1}{2}})\|^2+(\sigma^{\frac{1}{2}})^2~~\bigg)dxd\phi
        \leq Ck.
        \end{align*}
       Hence
       \begin{align*}
          &\mu^{(1),k}\{(x,\beta,\phi)~|~~\|\beta\|\geq M\}
          \leq e^{-\rho^{(1)}k}\frac{k}{2}(Ck)^{2\rho^{(1)}+1}+Ce^{-\rho^{(1)}k}\\
          &\leq C e^{-\rho^{(1)}k} (k^{2\rho^{(1)}+2}+1)\to 0,
       \end{align*}
         as $k \to \infty$.
        Therefore, passing to a subsequence if necessary, we  have a measure $\mu^{(1)}$ on $\mathbb{T}^n\times\mathbb{R}^n\times\mathbb{T}^1$ such that
      \begin{align*}
         \mu^{(1),k}\rightharpoonup\mu^{(1)}.
      \end{align*}
     Letting $\sigma^{(1)}$ be the projection of $\mu^{(1)}$ into $\mathbb{T}^{n}\times\mathbb{T}^{1}$, we have
      \begin{align*}
         \sigma^{(1),k}\rightharpoonup\sigma^{(1)}.
      \end{align*}

      \noindent\textbf{Step 4.}
     Further, we discuss the convergence of $\sigma^k_{x_ix_i}$.
      Let $\sigma^{(2),k}:=\sigma^k_{x_ix_i}=\sigma^{(1),k}_{x_i}
      =k[H(x,D_xu^k,\phi)]_{x_ix_i}\sigma^k+k^2[H(x,D_xu^k,\phi)]^2_{x_i}\sigma^k$,
      $\mu^{(2),k}$ lift measure $\sigma^{(2),k}$ to space $\mathbb{T}^{n}\times\mathbb{R}^{n}\times\mathbb{T}^{1}$.
      Then,  for sufficiently large $M$,
         \begin{align*}
          &\mu^{(2),k}\{(x,\beta,\phi)~|~~\|\beta\|\geq M\}\\
          &=e^{-\rho^{(2)}k}\int_{\mathbb{T}^{n}\times\mathbb{T}^{1}}
          (\sigma^k)^{\rho^{(2)}}\bigg(k[H(x,D_xu^k,\phi)]_{x_ix_i}\sigma^k
          +k^2[H(x,D_xu^k,\phi)]^2_{x_i}\sigma^k\bigg)dxd\phi\\
          &=e^{-\rho^{(2)}k}\int_{\mathbb{T}^{n}\times\mathbb{T}^{1}}
          -(\rho^{(2)}+1)k^2(\sigma^k)^{\rho^{(2)}+1}[H(x,D_xu^k,\phi)]^2_{x_i}
          +k^2(\sigma^k)^{\rho^{(2)}+1}[H(x,D_xu^k,\phi)]^2_{x_i}\sigma^k\bigg)dxd\phi\\
          &=-e^{-\rho^{(2)}k}\int_{\mathbb{T}^{n}\times\mathbb{T}^{1}}
          \rho^{(2)}k^2(\sigma^k)^{\rho^{(2)}+1}[H(x,D_xu^k,\phi)]^2_{x_i}dxd\phi,
         \end{align*}
         where $\rho^{(2)}$  is a finite number and we will provide the specific form below.
         On the one hand, because $[H(x,D_xu^k,\phi)]_{x_i}=H_{x_i}+H_{y_j}u^k_{x_jx_i}$ and
         $\max_{(x,\phi)\in\mathbb{T}^{n}\times\mathbb{T}^{m}}\|D_xu^k\|\leq C$,
        we have
         \begin{align*}
           |H_{x_i}|\leq C(1+\|D_xu^k\|^2)\leq C,  |H_{y_i}|\leq C(1+\|D_xu^k\|)\leq C.
         \end{align*}
         On the other hand, $u^k$ is the smooth solution of
         \begin{align*}
           {\rm div}_x[e^{k[H(x,D_xu^k,\phi)}D_yH(x,D_xu^k,\phi)]=0,
         \end{align*}
         that is,
         $(kH_{y_i}H_{y_j}+H_{y_iy_j})u^k_{x_jx_i}+(kH_{x_i}H_{y_i}+H_{y_ix_i})=0$.
         Based on  the Schauder estimate of elliptic equations, we obtain  the bounded estimate of $u^k_{x_ix_j}$, refer to \cite{G-T,L,X}.
         Then,
         \begin{align*}
          \mu^{(2),k}\{(x,\beta,\phi)~|~~\|\beta\|\geq M\}
           \leq Ce^{-\rho^{(2)}k}\rho^{(2)}k^2
           \int_{\mathbb{T}^{n}\times\mathbb{T}^{1}}
          (\sigma^k)^{\rho^{(2)}+1}dxd\phi.
         \end{align*}
         Based on Sobolev's embedded  inequality, that is, $||(\sigma^k)^{\frac{1}{2}}||_{L^{2(\rho^{(2)}+1)}(\mathbb{T}^{n}\times\mathbb{T}^{1})} \leq C ||(\sigma^k)^{\frac{1}{2}}||_{W^{1,2}(\mathbb{T}^{n}\times\mathbb{T}^{1})} \leq Ck$, and
        \begin{align*}
            \rho^{(2)}=\left\{
                \begin{array}{lr}
                  \frac{2}{n-1}, & n\geq2; \\
                  {\rm any ~positive ~number}, & n=1.
                \end{array}
             \right.
        \end{align*}
        Hence
        \begin{align*}
         \mu^{(2),k}\{(x,\beta,\phi)~|~~\|\beta\|\geq M\}
           \leq Ce^{-\rho^{(2)}k}\rho^{(2)}k^2(Ck)^{\rho^{(2)}+1}
           =Ce^{-\rho^{(2)}k}\rho^{(2)}k^{\rho^{(2)}+3} \to 0,
        \end{align*}
         as $k \to \infty$.
        Therefore, similarly, passing to a subsequence if necessary, we  have
         \begin{align*}
         \mu^{(2),k}\rightharpoonup\mu^{(2)},
      \end{align*}
      and
       \begin{align*}
         \sigma^{(2),k}\rightharpoonup\sigma^{(2)},
      \end{align*}
      where $\mu^{(2)}$  is a measure on $\mathbb{T}^n\times\mathbb{R}^n\times\mathbb{T}^1$, and   $\sigma^{(2)}$ be the projection of $\mu^{(2)}$ into $\mathbb{T}^{n}\times\mathbb{T}^{1}$.

      \noindent\textbf{Step 5.}
     We discuss the convergence of $\sigma^k_{x_ix_ix_i}$.
      Let
      \begin{align*}
          \sigma^{(3),k}:=\sigma^k_{x_ix_ix_i}=\sigma^{(2),k}_{x_i}
          &=k[H(x,D_xu^k,\phi)]_{x_ix_ix_i}\sigma^k
          +k^2[H(x,D_xu^k,\phi)]_{x_ix_i}[H(x,D_xu^k,\phi)]_{x_i}\sigma^k\\
          &+2k^2[H(x,D_xu^k,\phi)]_{x_i}[H(x,D_xu^k,\phi)]_{x_ix_i}\sigma^k
          +k^3[H(x,D_xu^k,\phi)]^3_{x_i}\sigma^k,
          \end{align*}
      $\mu^{(3),k}$ lift measure $\sigma^{(3),k}$ to space $\mathbb{T}^{n}\times\mathbb{R}^{n}\times\mathbb{T}^{1}$.
      Then,  for sufficiently large $M$,
         \begin{align*}
          &\mu^{(3),k}\{(x,\beta,\phi)~|~~\|\beta\|\geq M\}\\
          &=e^{-\rho^{(3)}k}\int_{\mathbb{T}^{n}\times\mathbb{T}^{1}}
          (\sigma^k)^{\rho^{(3)}+1}\bigg(k[H(x,D_xu^k,\phi)]_{x_ix_ix_i}
          +k^2[H(x,D_xu^k,\phi)]_{x_ix_i}[H(x,D_xu^k,\phi)]_{x_i}\bigg)dxd\phi\\
          &+e^{-\rho^{(3)}k}\int_{\mathbb{T}^{n}\times\mathbb{T}^{1}}
          (\sigma^k)^{\rho^{(3)}+1}\bigg(2k^2[H(x,D_xu^k,\phi)]_{x_i}[H(x,D_xu^k,\phi)]_{x_ix_i}
          +k^3[H(x,D_xu^k,\phi)]^3_{x_i}
          \bigg)dxd\phi\\
          &=-e^{-\rho^{(3)}k}\int_{\mathbb{T}^{n}\times\mathbb{T}^{1}}
          \rho^{(3)}k^2 (\sigma^k)^{\rho^{(3)}+1}[H(x,D_xu^k,\phi)]_{x_ix_i}[H(x,D_xu^k,\phi)]_{x_i}dxd\phi\\
          &=-e^{-\rho^{(3)}k}\int_{\mathbb{T}^{n}\times\mathbb{T}^{1}}
          \rho^{(3)}k^3(\sigma^k)^{\rho^{(3)}}[H(x,D_xu^k,\phi)]^3_{x_i}dxd\phi
         \end{align*}
         where $\rho^{(3)}$  is a finite number and we will provide the specific form below.
         because $[H(x,D_xu^k,\phi)]_{x_i}=H_{x_i}+H_{y_j}u^k_{x_jx_i}$,
         $[H(x,D_xu^k,\phi)]_{x_ix_i}=H_{x_ix_i}+H_{x_iy_j}u^k_{x_jx_i}+H_{x_iy_j}u^k_{x_jx_i}
         +H_{y_jy_l}u^k_{x_lx_i}u^k_{x_jx_i}+H_{y_j}u^k_{x_jx_ix_i}$,
         where $l:=1,2,\cdots,n$.
         According to the Schauder estimate of elliptic equations, we get  $u^k_{x_ix_j}$ is H\"{o}lder continuous, and then $u^k_{x_jx_ix_i}$ is bounded.
         Combining with the results of step 3, H\"{o}lder inequality and Sobolev embedded inequality, we obtain
        \begin{align*}
         \mu^{(3),k}\{(x,\beta,\phi)~|~~\|\beta\|\geq M\}
           \leq Ce^{-\rho^{(3)}k}\rho^{(3)}k^3(Ck)^{\rho^{(3)}+1}
           =Ce^{-\rho^{(3)}k}\rho^{(3)}k^{\rho^{(3)}+4} \to 0,
        \end{align*}
         as $k \to \infty$,
         where \begin{align*}
            \rho^{(3)}=\left\{
                \begin{array}{lr}
                  \frac{2}{n-1}, & n\geq2; \\
                  {\rm any ~positive ~number}, & n=1.
                \end{array}
             \right.
        \end{align*}
        Therefore, passing to a subsequence if necessary, we  have
         \begin{align*}
         \mu^{(3),k}\rightharpoonup\mu^{(3)},
      \end{align*}
      and
       \begin{align*}
         \sigma^{(3),k}\rightharpoonup\sigma^{(3)},
      \end{align*}
      where $\mu^{(3)}$  is a measure on $\mathbb{T}^n\times\mathbb{R}^n\times\mathbb{T}^1$, and   $\sigma^{(3)}$ be the projection of $\mu^{(3)}$ into $\mathbb{T}^{n}\times\mathbb{T}^{1}$.

         \noindent\textbf{Step 6.}
          According to the above discussion,  the regularity theory of elliptic equations and Vitali's convergence theorem, we can conclude that the limits
         \begin{align*}
          \lim_{k\to\infty}\int_{\mathbb{T}^n}G(z)N_{x_i}u^k(x,z)dz
          =\lim_{k\to\infty}u^k_{x_i}(x),
        \end{align*}
        and
         \begin{align*}
         \lim_{k\to\infty}\int_{\mathbb{T}^n}G(z)N_{x_ix_i}u^k(x,z)dz
         =\lim_{k\to\infty} u^k_{x_ix_i}(x)
        \end{align*}
         exist.
         Therefore,  the weak KAM solution $u(x,P,\phi)$ is also the classical solution of equation
         \begin{align*}
          \sigma u_{x_ix_i}+\sigma_{x_i} u_{x_i}=0.
         \end{align*}

            \end{proof}
          \section{Convergence}\label{sec:3}
       We obtain some convergence through a bunch of techniques in this section, and then get some limits.
       These limits  may not be better than the property of  original sequence in some respects, but they have special  advantages  to complete our conclusion.
       \begin{lemm}\label{lem-3.1}
				For every $\varrho\in\mathbb{N}^+$, there is a constant $C:=C(P,\varrho)$ such that
                       \begin{align*}
                          \begin{array}{lll}
                               \displaystyle  \mathop {\sup}_{k} \|u^k\|_{L^\varrho(\mathbb{T}^{n}\times\mathbb{T}^{m})}\leq C,
                             & \displaystyle  \mathop {\sup}_{k} \|D_xu^k\|_{L^\varrho(\mathbb{T}^{n}\times\mathbb{T}^{m})} \leq C,
                          \end{array}
                       \end{align*}
                where $\|u^k\|_{L^\varrho(\mathbb{T}^{n}\times\mathbb{T}^{m})}=\bigg(\int_{\mathbb{T}^{n}\times\mathbb{T}^{m}} |u^k|^\varrho dxd\phi\bigg)^{1/\varrho}$, $\|D_xu^k\|_{L^\varrho(\mathbb{T}^{n}\times\mathbb{T}^{m})}=\bigg(\int_{\mathbb{T}^{n}\times\mathbb{T}^{m}} \|D_xu^k\|^\varrho dxd\phi\bigg)^{1/\varrho}$.
		\end{lemm}
        \begin{proof}
        The proofs can be referred to \cite{E03,N-J-L}.
        \end{proof}

        According to  Lemma \ref{lem-3.1},  if necessary passing to a subsequence, there exists limit $u(x,P,\phi)$ such that the sequence $\{u^k\}_{k=1}^{\infty}$ converges uniformly to $u(x,P,\phi)$ on $\mathbb{T}^{n}\times\mathbb{T}^{m}$, denoted by
       \begin{align*}
          u^k(x,P,\phi)\rightarrow u(x,P,\phi):=x\cdot P+v(x,P,\phi)
       \end{align*}
        and for each $\varrho\in\mathbb{N}^+$, the sequence $\{D_xu^k\}_{k=1}^{\infty}$ weakly converges to $D_xu$ in $L^\varrho(\mathbb{T}^{n}\times\mathbb{T}^{m})$, denoted by
       \begin{align*}
          D_xu^k(x,P,\phi)\rightharpoonup D_xu(x,P,\phi):=P+D_xv(x,P,\phi).
       \end{align*}
        Next, we will talk about another kind of convergence.
        Firstly, we give  definition of the tightness, and then  prove that the probability  measures sequence $\{\mu^{k}\}^{\infty}_{k=1}$ is tight.

        \begin{definition}
          Let $(\Omega, T)$ be a Hausdorff space, and let $\Sigma$  be a $\sigma$-algebra on $\Omega$ that contains the topology $T$. Let $M$ be a collection of measures defined on $\Sigma$. The collection $M$ is called tight if, for any $\varepsilon >0$, there is a compact subset $K_{\varepsilon }$ of $\Omega$ such that, for all measures $\eta \in M$,
          \begin{align*}
             \eta (\Omega\setminus K_{\varepsilon }) < \varepsilon
          \end{align*}
          hold.
       \end{definition}

      \begin{prop}\label{pro-3.2}
         The probability measure sequence $\{\mu^{k}\}^{\infty}_{k=1}$ defined above on $\mathbb{T}^n\times\mathbb{R}^n\times\mathbb{T}^m$  is tight.
      \end{prop}

       \begin{proof}

       Let us leave out the superscript $k$ for now. On the basis of   ${\rm div}_x[\sigma^kD_yH(x,P+D_xv^k,\phi)]=0$, we can obtain
             \begin{align}\label{2-1}
                 0&=\int_{\mathbb{T}^{n}\times\mathbb{T}^{m}} (\sigma H_{y_i})_{x_i}u_{x_ix_j}dxd\phi\nonumber\\
                  &=\int_{\mathbb{T}^{n}\times\mathbb{T}^{m}}  (H_{y_ix_i}+H_{y_iy_j}u_{x_jx_i})\sigma u_{x_ix_j}+\sigma k [H(x,D_xu,\phi)]_{x_i}H_{y_i}u_{x_ix_j}dxd\phi.
              \end{align}
       Obviously, we have the following
             \begin{align*}
                \begin{array}{cc}
                 (\sigma)_{x_i}=k \sigma[H(x,D_xu,\phi)]_{x_i},
                 & [H(x,D_xu,\phi)]_{x_i}=H_{x_i}+H_{y_j}u_{x_jx_i}.
                \end{array}
             \end{align*}
       Hence we have
            \begin{align}\label{2-2}
               &\int_{\mathbb{T}^{n}\times\mathbb{T}^{m}} \sigma k [H(x,D_xu,\phi)]_{x_i} H_{y_i}u_{x_ix_j}dxd\phi
               =\int_{\mathbb{T}^{n}\times\mathbb{T}^{m}}  k\sigma [H(x,D_xu,\phi)]_{x_i} \bigg([H(x,D_xu,\phi)]_{x_j}-H_{x_j}\bigg)dxd\phi\nonumber\\
               &=\int_{\mathbb{T}^{n}\times\mathbb{T}^{m}}\sigma k [H(x,D_xu,\phi)]_{x_i} [H(x,D_xu,\phi)]_{x_j}dxd\phi
               +\int_{\mathbb{T}^{n}\times\mathbb{T}^{m}} \sigma(H_{x_jy_i}u_{x_ix_j}+H_{x_jx_i} )dxd\phi\nonumber\\
               &=\int_{\mathbb{T}^{n}\times\mathbb{T}^{m}}k [H(x,D_xu,\phi)]_{x_i} [H(x,D_xu,\phi)]_{x_j}d\sigma
               +\int_{\mathbb{T}^{n}\times\mathbb{T}^{m}} (H_{x_jy_i}u_{x_ix_j}+H_{x_jx_i}) d\sigma.
       \end{align}
       We combine (\ref{2-1}) and (\ref{2-2}) to  yield
        \begin{align}\label{2-4}
        \int_{\mathbb{T}^{n}\times\mathbb{T}^{m}} k [H(x,D_xu,\phi)]_{x_i} [H(x,D_xu,\phi)]_{x_j}d\sigma
        =\int_{\mathbb{T}^{n}\times\mathbb{T}^{m}} -H_{y_iy_j}u_{x_jx_i}u_{x_ix_j}-H_{y_ix_i}u_{x_ix_j}-H_{x_jy_i}u_{x_ix_j}-H_{x_jx_i}d\sigma.
        \end{align}
        Because of  Lemma \ref{lem-3.1}, (\ref{2-4}) and hypothesis of $H$,
        \begin{align*}
        \int_{\mathbb{T}^{n}\times\mathbb{T}^{m}} k [H(x,D_xu,\phi)]_{x_i}^2d\sigma
        \leq  -\frac{\gamma}{2}\int_{\mathbb{T}^{n}\times\mathbb{T}^{m}}|u_{x_ix_j}|^2d\sigma+C\bigg(\int_{\mathbb{T}^{n}\times\mathbb{T}^{m}}|u_{x_ix_j}|^2d\sigma\bigg)^{1/2}+C\leq C.
        \end{align*}
        Then
        \begin{align*}
        \int_{\mathbb{T}^{n}\times\mathbb{T}^{m}}\bigg(\|D_x(\sigma^{\frac{1}{2}})\|^2+(\sigma^{\frac{1}{2}})^2~~\bigg)dxd\phi
        = \frac{k^2}{4}\int_{\mathbb{T}^{n}\times\mathbb{T}^{m}}[H(x,D_xu,\phi)]_{x_i}^2d\sigma +1
        \leq Ck.
        \end{align*}
           Further,
         \begin{align}\label{2-A}
           \bigg(\int_{\mathbb{T}^{n}\times\mathbb{T}^{m}} \sigma^{~\rho} d\sigma\bigg)^{1/(1+\rho)}
           \leq C \int_{\mathbb{T}^{n}\times\mathbb{T}^{m}}\bigg(~~\|D_x(\sigma^{\frac{1}{2}})\|^2+(\sigma^{\frac{1}{2}})^2~~\bigg)dxd\phi
           \leq Ck,
         \end{align}
        where the first inequality of the above formula is due to Sobolev's embedded  inequality, that is, $||\sigma^{\frac{1}{2}}||_{L^{2(1+\rho)}(\mathbb{T}^{n}\times\mathbb{T}^{m})} \leq C ||\sigma^{\frac{1}{2}}||_{W^{1,2}(\mathbb{T}^{n}\times\mathbb{T}^{m})}$, and
        \begin{align*}
            \rho=\left\{
                \begin{array}{lr}
                  \frac{2}{n+m-2}, & n+m\geq3; \\
                  {\rm any ~positive ~number}, & n+m=2.
                \end{array}
             \right.
        \end{align*}
        Finally,  we restore the superscript $k$, for sufficiently large $M$,
         \begin{align*}
          \mu^{k}\{(x,\beta,\phi)~|~~\|\beta\|\geq M\} &\leq \sigma^{k}\{(x,\phi)~|~H(x,D_xu^k,\phi)\geq \bar{H}^k(P)+1\}
          \leq\int_{\mathbb{T}^{n}\times\mathbb{T}^{m}} e^{~\rho k(H(x,D_xu^k,\phi)-\bar{H}^k(P)-1)}d\sigma^k\\
          &= e^{-\rho k}\int_{\mathbb{T}^{n}\times\mathbb{T}^{m}} (\sigma^{k})^{\rho} d\sigma^{k} \leq C e^{-k\rho}k^{1+\rho} \to 0,
         \end{align*}
         as $k \to \infty$.
    \end{proof}

      Passing to a subsequence if necessary, we  have a probability measure $\mu$ on $\mathbb{T}^n\times\mathbb{R}^n\times\mathbb{T}^m$ such that
      \begin{align*}
         \mu^k\rightharpoonup\mu.
      \end{align*}
      Similarly, letting $\sigma$ be the projection of $\mu$ into $\mathbb{T}^{n}\times\mathbb{T}^{m}$, we have
      \begin{align*}
         \sigma^k\rightharpoonup\sigma.
      \end{align*}
    \section{Properties of limits}\label{sec:4}
       In the last section, we obtained the limit of function $u$ and the probability measure  $\mu$.
       Next, we will discuss some  properties about $u$ and $\mu$.
       To  prove these properties better, we firstly present some preparatory lemmas.
    \begin{lemm}\label{lem-4.1}
            Suppose $H>0$, then the following assertions are true:
                   \begin{description}
                         \item (i)~~  $H(x,D_x u,\phi)\leq\bar{H}(P)$ almost everywhere in $\mathbb{T}^{n}\times\mathbb{T}^{m}$.
                         Further, we obtain
                         \begin{align*}
                            H(x,D_xu^\varepsilon,\phi) \leq \bar{H}(P)+O(\varepsilon)-\chi_\varepsilon(x,\phi),
                         \end{align*}
                       where $u^\varepsilon:=\rho_\varepsilon\ast u$, $\rho_\varepsilon$  is the standard mollified kernel, $\chi_\varepsilon(x,\phi)=\int_{\mathbb{R}^{n}\times\mathbb{R}^{m}}\rho_\varepsilon\|D_xu(\tilde{x},\tilde{\phi})-D_xu^\varepsilon\|^2d\tilde{x}d\tilde{\phi}$;
                         \item (ii)~~$\bar{H}(P)=\displaystyle\mathop{\lim}_{k\to\infty}\bar{H}^k(P)
                             =\displaystyle\mathop{\lim}_{k\to\infty}\frac{1}{k}\log\bigg(\int_{\mathbb{T}^{n}\times\mathbb{T}^{m}}
                                      e^{kH(x,D_xu^k,\phi)}dxd\phi\bigg)$.
                   \end{description}
          \end{lemm}
          A proof of the Lemma \ref{lem-4.1} can be found in the Appendix.
          \begin{lemm}\label{lem-4.2}
          Inf-max formula holds. That is,
          \begin{align*}
                \bar{H}(P)=\mathop {\inf}_{\hat{v}\in C^1(\mathbb{T}^{n})\times C(\mathbb{T}^{m})}\mathop {\max}_{(x,\phi)\in\mathbb{T}^{n}\times\mathbb{T}^{m}}H(x,P+D_x\hat{v},\phi).
          \end{align*}
          \end{lemm}

          See \cite{E03,N-J-L} for proof of  the Lemma \ref{lem-4.2}.

          \begin{lemm}\label{lem-4.3}
          The measure $\mu$ is closed, that is,
           \begin{align*}
           \int_{\mathbb{T}^{n}\times\mathbb{R}^n\times\mathbb{T}^{m}}\beta\cdot D_x\varphi d\mu=0
           \end{align*}
            for all $\varphi\in C^1(\mathbb{T}^{n}\times\mathbb{T}^{m})$.
          \end{lemm}
          \begin{proof}
           \begin{align*}
            \int_{\mathbb{T}^{n}\times\mathbb{R}^n\times\mathbb{T}^{m}}\beta\cdot D_x\varphi d\mu
              &=\displaystyle\mathop{\lim}_{k\to\infty}\int_{\mathbb{T}^{n}\times\mathbb{T}^{m}} D_yH(x,D_xu^k,\phi)\cdot D_x\varphi d\sigma^k\\
              &=\displaystyle\mathop{\lim}_{k\to\infty}\int_{\mathbb{T}^{n}\times\mathbb{T}^{m}} D_yH(x,D_xu^k,\phi)\cdot \sigma^k d\varphi d\phi\\
               &=\displaystyle\mathop{\lim}_{k\to\infty}-\int_{\mathbb{T}^{n}\times\mathbb{T}^{m}} {\rm div}_x[\sigma^kD_yH(x,D_xu^k,\phi)]~\varphi~ dxd\phi=0.
             \end{align*}
          \end{proof}
         After all the preparations are made, we prove that $\mu$ is a minimal measure (Proposition \ref{pro-4.4}).
         Additionally,  we demonstrate that $u$ satisfies  a certain smoothness (Lemma \ref{lem-4.5} and Proposition \ref{pro-4.7}) and a Hamilton-Jacobi equation in the sense of the minimal measure (Theorem \ref{th3}).
           \begin{prop}\label{pro-4.4}
         The limit
        \begin{align*}
            \displaystyle\mathop{\lim}_{k\to\infty} \int_{\mathbb{T}^{n}\times\mathbb{T}^{m}} H(x,D_xu^k,\phi) d\sigma^k=\bar{H}(P)
        \end{align*}
        holds. The measure $\mu$ is a minimal measure, and
        \begin{align*}
            \int_{\mathbb{T}^{n}\times\mathbb{R}^n \times\mathbb{T}^{m}} L(x,\beta,\phi)d\mu =\bar{L}(Q),
       \end{align*}
       where $Q=\int_{\mathbb{T}^{n}\times\mathbb{R}^n\times\mathbb{T}^{m}} \beta~ d\mu$, $\beta=D_yH(x,D_xu^k,\phi)$.
       \end{prop}

          \begin{proof}
         First, we claim that
         \begin{align}\label{3-1}
            \bar{H}(P)
            \leq \displaystyle\mathop{\liminf}_{k\to\infty} \int_{\mathbb{T}^{n}\times\mathbb{T}^{m}} H(x,D_xu^k,\phi)d\sigma^k.
         \end{align}
         Let $\Pi:=\{(x,\phi)~|~H(x,D_xu^k,\phi) \geq \bar{H}^k(P)-\lambda~{\rm and}~\lambda>0\}$.
         We compute for each $\lambda$ that
        \begin{align*}
           \bar{H}^k(P)&= \int_{\mathbb{T}^{n}\times\mathbb{T}^{m}} \bar{H}^k(P) d\sigma^k
                            = \int_{\mathbb{T}^{n}\times\mathbb{T}^{m}\cap\Pi} \bar{H}^k(P) d\sigma^k+\int_{\mathbb{T}^{n}\times\mathbb{T}^{m}\backslash\Pi} \bar{H}^k(P) d\sigma^k\\
            &\leq\int_{\mathbb{T}^{n}\times\mathbb{T}^{m}\cap\Pi}   H(x,D_xu^k,\phi)d\sigma^k+\lambda+\int_{\mathbb{T}^{n}\times\mathbb{T}^{m}\backslash\Pi}\bar{H}^k(P)e^{-k\lambda} dxd\phi\\
            &\leq \int_{\mathbb{T}^{n}\times\mathbb{T}^{m}}  H(x,D_xu^k,\phi)d\sigma^k+\lambda+Ce^{-k\lambda}.
       \end{align*}
       Because of Lemma \ref{lem-4.1} (i),
      \begin{align*}
          \bar{H}(P)=\displaystyle\mathop{\lim}_{k\to\infty}\bar{H}^k(P)
          \leq \displaystyle\mathop{\liminf}_{k\to\infty} \int_{\mathbb{T}^{n}\times\mathbb{T}^{m}} H(x,D_xu^k,\phi) d\sigma^k+\lambda,
     \end{align*}
     and  since $\lambda>0$ is arbitrary, (\ref{3-1}) follows.
     Then,  we apply (\ref{3-1}) to calculate
      \begin{align}\label{3-2}
          \int_{\mathbb{R}^{m}\times\mathbb{T}^{n}\times\mathbb{T}^{m}} L(x,\beta,\phi)d\mu
          &=\displaystyle\mathop{\lim}_{k\to\infty}\int_{\mathbb{T}^{n}\times\mathbb{T}^{m}} L(x,D_y(x,D_xu^k,\phi),\phi)d\sigma^k\nonumber\\
          &=\displaystyle\mathop{\lim}_{k\to\infty}\int_{\mathbb{T}^{n}\times\mathbb{T}^{m}} D_xu^k\cdot D_yH(x,D_xu^k,\phi)-H(x,D_xu^k,\phi)d\sigma^k\nonumber\\
          &=\displaystyle\mathop{\lim}_{k\to\infty} \int_{\mathbb{T}^{n}\times\mathbb{T}^{m}} P\cdot D_yH(x,D_xu^k,\phi) -H(x,D_xu^k,\phi)d\sigma^k.\nonumber\\
          &\leq P\cdot Q-\bar{H}(P)\leq\displaystyle\mathop{\sup}_{R}\{R\cdot Q-\bar{H}(R)\}=\bar{L}(Q).
     \end{align}
     But Mather's minimization principle asserts that
     \begin{align*}
         \bar{L}(Q)=\inf \bigg\{ \int_{\mathbb{R}^m\times\mathbb{T}^{n}\times\mathbb{T}^{m}} L(x,\beta,\phi)d\mu~|~\mu~{\rm is~flow~invariant},~Q=\int_{\mathbb{R}^n\times\mathbb{T}^{n}\times\mathbb{T}^{m}} \beta~ d\mu~\bigg\}.
     \end{align*}
     Hence the inequality in (\ref{3-2})  becomes  equality, and  Proposition \ref{pro-4.4} is completely proven.
     \end{proof}

     Before  giving Theorem \ref{th3}, we give the  smoothness of $u(x,P,\phi)$ with respect to $x$.
      \begin{lemm}\label{lem-4.5}
      The limit of  function $u(x,P,\phi):=x\cdot P+v(x,P,\phi)$ is  differentiable $\sigma$-almost everywhere with respect to $x$, where $\sigma$ is a minimal measure.
      \end{lemm}
The detailed proof can be found in the Appendix.
        \begin{theo}\label{th3}
        The  function $u(x,P,\phi):=x\cdot P+v(x,P,\phi)$ satisfies
         \begin{align}\label{3-C}
            H(x,D_xu,\phi)=\bar{H}(P)
        \end{align}
        $\sigma-$  almost everywhere,
         \begin{align}\label{3-D}
            {\rm div}_x[\sigma D_yH(x,D_xu,\phi)]=0
         \end{align}
         in  $\mathbb{T}^{n}\times\mathbb{T}^{m}$, where $\sigma$ is a minimal measure.
         \end{theo}

         \begin{proof}

       We use proof by contradiction to obtain (\ref{3-C}).
       Assume $H(x,D_xu,\phi)<\bar{H}(P)$ on a set of minimal measure  $\sigma$,
       we set now $r=D_yH(x,D_xu,\phi)$, then
         \begin{align*}
         0&=\frac{\gamma}{2}\int_{\mathbb{R}^{m}\times\mathbb{T}^{n}\times\mathbb{T}^{m}}\|\beta-D_yH(x,D_xu,\phi)\|^2d\mu\nonumber\\
         &\leq \int_{\mathbb{R}^{m}\times\mathbb{T}^{n}\times\mathbb{T}^{m}} L(x,\beta,\phi)-L(x,D_yH(x,D_xu,\phi),\phi)-D_\beta L(x,D_yH(x,D_xu,\phi),\phi)\cdot(\beta-D_yH(x,D_xu,\phi))d\mu\nonumber\\
         &=\int_{\mathbb{R}^{m}\times\mathbb{T}^{n}\times\mathbb{T}^{m}} L(x,\beta,\phi) + H(x,D_xu,\phi)-D_xu\cdot D_yH(x,D_xu,\phi)-D_xu\cdot (\beta-D_yH(x,D_xu,\phi)) d\mu\nonumber\\
         &=\int_{\mathbb{R}^{m}\times\mathbb{T}^{n}\times\mathbb{T}^{m}} L(x,\beta,\phi) + H(x,D_xu,\phi)-D_xu \cdot \beta d\mu\\
         &<\bar{L}(Q)+\bar{H}(P)-P\cdot Q=0.
      \end{align*}
       This is contradictory, so
       \begin{align*}
           H(x,D_xu,\phi)\geq\bar{H}(P)
       \end{align*}
      on a set of positive $\sigma$ measure and in light of Lemma \ref{lem-4.1} (ii), we can obtain
       \begin{align*}
           H(x,D_xu,\phi)=\bar{H}(P)
       \end{align*}
        $\sigma$-almost everywhere. For each $\varphi\in C^1(\mathbb{T}^{n}\times\mathbb{T}^{m})$, due to Lemma \ref{lem-4.3} ($\mu$ is closed measure) and
        $\beta=D_yH(x,D_xu,\phi)$ in $\mu$-almost everywhere, we have
        \begin{align*}
            0&=\int_{\mathbb{R}^{m}\times\mathbb{T}^{n}\times\mathbb{T}^{m}} \beta\cdot D_x\varphi d\mu
             =\int_{\mathbb{T}^{n}\times\mathbb{T}^{m}}  D_yH(x,D_xu,\phi)\cdot D_x\varphi d\sigma
             =\int_{\mathbb{T}^{n}\times\mathbb{T}^{m}}  D_yH(x,D_xu,\phi)\cdot \sigma d\varphi d\phi\\
             &=-\int_{\mathbb{T}^{n}\times\mathbb{T}^{m}} {\rm div}_x[\sigma D_yH(x,D_xu,\phi)] ~\varphi~ dxd\phi.
        \end{align*}
         Hence, (\ref{3-D}) holds.
       \end{proof}
         Finally, according to the continuous dependence of the solution on the parameters, we obtain that the weak KAM solution $u(x,P,\phi)$ is continuous with respect to $\phi$.
          That is,
        \begin{prop}\label{pro-4.7}
          The weak KAM solution $u(x,P,\phi)$    is continuous with respect to $\phi$.
        \end{prop}

        \begin{proof}
         First, for each $k\in\mathbb{N}^+$, we claim that  $v^k_\phi$  is continuous with respect to $\phi$.
         Fixed $\phi_0\geq0$, we assume that positive sequence $\{\phi_i\}\rightarrow \phi_0$, such that $\displaystyle \mathop{\lim}_{\phi_i \to \phi_0} v^k_{\phi_i}=\hat v^k_{\phi_0}$.
         To get the continuity of $v^k_{\phi}$ with respect to $\phi$, it suffices to prove $v^k_{\phi_0}=\hat v^k_{\phi_0}$.
         Otherwise, we have
         \begin{align}\label{l4-1}
			\|\hat v^k_{\phi_0}-v^k_{\phi_0}\|_{L^2(\mathbb{T}^{n}\times\mathbb{T}^{m})}
            +\|D_x\hat v^k_{\phi_0}-D_xv^k_{\phi_0}\|_{L^2(\mathbb{T}^{n}\times\mathbb{T}^{m})}\geq a
	     \end{align}
         for some $a>0$. For fixed $k$, $v^k_{\phi_0}$ and $\hat v^k_{\phi_0}$ are  minimizers of $I_k[\cdot]$ as $\phi=\phi_0$.
		 Let $g(t)=I_k[(1-t)\hat v^k_{\phi_0}+tv^k_{\phi_0}]$, obviously, $g(0)=I_k[\hat v^k_{\phi_0}]=I_k[v^k_{\phi_0}]=g(1)$.
			By Taylor's formula, $g(1)$ has the following form
			\begin{align}\label{l4-2}
				g(1)=g(0)+g'(0)+\frac{g''(\xi)}{2},
			\end{align}
               where $\xi\in(0,1)$.
			Since $\hat v^k_{\phi_0}$ is the minimizer of $I_k[\cdot]$,  $g'(0)=\delta I_k[(1-t)\hat v^k_{\phi_0}+tv^k_{\phi_0}]\textbf{~\big|}_{t=0}=0$, which together with (\ref{l4-2}) implies $g''(\xi)=0$.
			Let $\tilde{s}=(1-\xi)\hat v^k_{\phi_0}+\xi v^k_{\phi_0}$, $s=v^k_{\phi_0}-\hat v^k_{\phi_0}$,
			then
			\begin{align*}
				g''(\xi)&=\delta^2I_k[(1-t)\hat v^k_{\phi_0}+tv^k_{\phi_0}]\textbf{~\big|}_{t=\xi}\\
				&=\int_{\mathbb{T}^{n}\times\mathbb{T}^{m}}k^2 e^{kH(x,D_x\tilde{s},\phi)}
                (H_{y_i}s_{x_i})^2
                +ke^{kH(x,D_x\tilde{s},\phi)}(H_{y_iy_j}s_{x_i}s_{x_j})dxd\phi
			\end{align*}
			and $H$ is evaluated at $(x,D_x\tilde{s},\phi)$. Apparently,
            \begin{align*}
            \int_{\mathbb{T}^{n}\times\mathbb{T}^{m}}k^2 e^{kH(x,D_x\tilde{s},\phi)}
                (H_{y_i}s_{x_i})^2dxd\phi \geq0.
           \end{align*}
			According to (H1), we have
           \begin{align*}
            \int_{\mathbb{T}^{n}\times\mathbb{T}^{m}}ke^{kH(x,D_x\tilde{s},\phi)}(H_{y_iy_j}s_{x_i}s_{x_j})dxd\phi
            \geq\frac{\lambda}{2} \int_{\mathbb{T}^{n}\times\mathbb{T}^{m}}ke^{kH}|D_xs|^2dxd\phi,
            \end{align*}
            and then, due to  $g''(\xi)=0$,
            \begin{align}\label{l4-3}
            \|D_xs\|_{L^2(\mathbb{T}^{n}\times\mathbb{T}^{m})}^2=0.
             \end{align}
            Next, by Poincar\'{e} inequality, we have
             \begin{align*}
             \|s\|_{L^2(\mathbb{T}^{n}\times\mathbb{T}^{m})}
             \leq C \|D_xs\|_{L^2(\mathbb{T}^{n}\times\mathbb{T}^{m})}
             \end{align*}
             for some constant $C$. Then, combining with (\ref{l4-3}), we obtain
             \begin{align*}
             0=s=v^k_{\phi_0}-\hat v^k_{\phi_0}.
             \end{align*}
            This is in contradiction with (\ref{l4-1}).
			Therefore, the weak KAM solution  $u(x,y,\phi)$  is  continuous with respect to $\phi$, as $k\rightarrow\infty$.
           \end{proof}
       Moreover, the weak KAM solution  $u(x,y,\phi)$ satisfies Aronsson's equation
       \begin{align*}
              H_{y_i}(x,D_xu,\phi)H_{y_j}(x,D_xu,\phi)u_{x_ix_j}+H_{x_i}(x,D_xu,\phi) H_{y_i}(x,D_xu,\phi)=0
       \end{align*}
       in the sense of viscosity solutions.
       We can refer to \cite{E03} for a detailed proof.

      \section{ Appendix \bf}\label{sec:6}

       \subsection{\textbf{Existence and uniqueness of the solution}}
     We discuss  the smooth solution of
     \begin{align}\label{4-1}
        {\rm div}_x[e^{k[H(x,P+D_xv^k,\phi)}D_yH(x,P+D_xv^k,\phi)]=0
     \end{align}
     by the continuation method.
     \begin{prop} \label{prop-5.1}
     These exists  a unique  smooth solution to  the Euler-Lagrange equation (\ref{4-1}).
     \end{prop}
     \begin{proof}
     View $\phi$ as a parameter, and we consider the family of Hamilton functions
     \begin{align*}
        H_\tau(x,y,\phi):=\tau H(x,y,\phi)+(1-\tau)\frac{|y|^2}{2},
     \end{align*}
     and for each $k$ introduce the PDE
     \begin{align}\label{4-2}
         {\rm div}_x[e^{k[H_\tau(x,D_xu^k,\phi)}D_yH_\tau(x,D_xu^k,\phi)]=0,
     \end{align}
     where $u^k=P\cdot x+v^k(x,y,\phi)$, $v^k$ being $\mathbb{T}^n$-periodic.
     Define
     \begin{align*}
        \Gamma:=\{\tau ~|~\tau\in[0,1], (\ref{4-2}){\rm ~has~ a~ smooth~ solution}\}.
     \end{align*}
     Clearly, we have $0\in\Gamma$, corresponding to the solution $u^k=P\cdot x$, that is $v^k=0$.

       We claim that $\Gamma$ is closed. Since the following Lemma \ref{lem-5.2}  and $H_\tau$  satisfies the basic assumptions.
     We can obtain, for any $\tau\in\Gamma$,  the boundedness of the first derivative of the solution to (\ref{4-2}).
     Elliptic regularity theory implies that the solution is uniformly bounded in $\tau\in\Gamma$ for any order of derivatives.
     Thus, any convergent sequence has a subsequence whose corresponding solution of   (\ref{4-2}) converge uniformly, along with all derivatives.

       We claim that $\Gamma$ is open. The linearization of (\ref{4-2}) about $u^k$ is
       \begin{align*}
          Lw~:=-{\rm div}_x(e^{kH_\tau(x,D_xu^k,\phi)}(H_{\tau,y_iy_j}(x,D_xu^k,\phi)
              +kH_{\tau,y_i}(x,D_xu^k,\phi)H_{\tau,y_j}(x,D_xu^k,\phi))w_{x_j})_{x_i}.
       \end{align*}

     This is a symmetric, uniformly elliptic operator and null space consists of the constant. Using the Implicit function theorem, this yields a unique solution for nearby values of $\tau$.
      Because $u$ is continuous with respect to $\phi$ (Proposition \ref{pro-4.7}) and $\mathbb{T}^m$  is a compact manifold without boundary,
      ones can obtain $u$ is $\mathbb{T}^m$-periodic.

       Since $\Gamma$ is nonempty, closed and open, we have $\Gamma=[0,1]$ and in particular $1 \in \Gamma$. Thus (\ref{4-1}) has a smooth solution $u^k=x\cdot P+v^k(x,P,\phi)$.
       And further, if need to obtain  uniqueness of smooth solution, we require the condition $\int_{\mathbb{T}^{n}\times\mathbb{T}^{m}}v^k dxd\phi=0$. Because any such solution gives the unique minimizer of $I_k[\cdot]$.
     \end{proof}
     \begin{lemm} \label{lem-5.2}
     These exists a constant $C>0$ such that
     \begin{align*}
        \max_{(x,\phi)\in\mathbb{T}^{n}\times\mathbb{T}^{m}}\|D_xu^k\|\leq C
     \end{align*}
     for $k\in \mathbb{N}^+$.
     \end{lemm}
     \begin{proof}
        We drop the superscripts $k$, according to ${\rm div}_x[\sigma^kD_yH(x,P+D_xv^k,\phi)]=0$,
        \begin{align*}
           0&=\int_{\mathbb{T}^{n}\times\mathbb{T}^{m}}(\sigma H_{y_i})_{x_i}(\sigma^pH_{y_i})_{x_i}dxd\phi
            =\int_{\mathbb{T}^{n}\times\mathbb{T}^{m}}  ((\sigma_{x_i}H_{y_i}+\sigma(H_{y_i})_{x_i})(p\sigma^{p-1}\sigma_{x_i}H_{y_i}+\sigma^p(H_{y_i})_{x_i})dxd\phi\\
            &=\int_{\mathbb{T}^{n}\times\mathbb{T}^{m}} p\sigma^{p-1}(\sigma_{x_i}H_{y_i})^2+\sigma^p\sigma_{x_i}H_{y_i}(H_{y_i})_{x_i}
            +p\sigma^p(H_{y_i})_{x_i}\sigma_{x_i}H_{y_i}+\sigma^{p+1}((H_{y_i})_{x_i})^2dxd\phi\\
            &=\int_{\mathbb{T}^{n}\times\mathbb{T}^{m}} p\sigma^{p-1}(\sigma_{x_i}H_{y_i})^2+\sigma^{p+1}((H_{y_i})_{x_i})^2+(p+1)\sigma^p(H_{y_i})_{x_i}\sigma_{x_i}H_{y_i}dxd\phi\\
            &=:\int_{\mathbb{T}^{n}\times\mathbb{T}^{m}} A+B+Cdxd\phi.
        \end{align*}
        Clearly, $A$,$B\geq0$, and we further calculate that
        \begin{align*}
           C&=(p+1)\sigma^p(H_{y_i})_{x_i}\sigma_{x_i}H_{y_i}
            =(p+1)\sigma^p\sigma_{x_i}H_{y_i}(H_{y_iy_j}u_{x_jx_i}+H_{y_ix_i})\\
            &=\frac{p+1}{k} \sigma^{p-1}H_{y_iy_j}\sigma_{x_i}\sigma_{x_j}+(p+1)\sigma^p\sigma_{x_i}(H_{y_i}H_{y_ix_i}-H_{y_iy_j}H_{x_j})\\
            &\geq\frac{p+1}{k} \sigma^{p-1}\|D_x\sigma\|^2-C(p+1)\sigma^p\|D_x\sigma\|(1+\|D_xu\|^2),
        \end{align*}
        and so
        \begin{align*}
           \int_{\mathbb{T}^{n}\times\mathbb{T}^{m}} \sigma^{p-1}\|D_x\sigma\|^2dxd\phi\leq Ck\int_{\mathbb{T}^{n}\times\mathbb{T}^{m}} \sigma^p\|D_x\sigma\|(1+\|D_xu\|^2)dxd\phi.
        \end{align*}
        Then, by H\"{o}lder inequality,
        \begin{align*}
           \int_{\mathbb{T}^{n}\times\mathbb{T}^{m}} \sigma^{p-1}\|D_x\sigma\|^2dxd\phi\leq Ck^2\int_{\mathbb{T}^{n}\times\mathbb{T}^{m}} \sigma^{p+1}(1+\|D_xu\|^4)dxd\phi,
        \end{align*}
        Sobolev's inequality provides the estimate
        \begin{align*}
           \bigg(\int_{\mathbb{T}^{n}\times\mathbb{T}^{m}} \sigma^{(p+1)(1+\theta)}\bigg)^{1/{(1+\theta)}}
           &\leq C\int_{\mathbb{T}^{n}\times\mathbb{T}^{m}} \|D_x\sigma^{(p+1)/2}\|^2 +\sigma^{p+1}dxd\phi\\
           &=C(p+1)^2\int_{\mathbb{T}^{n}\times\mathbb{T}^{m}}\sigma^{p-1}\|D_x\sigma\|^2dxd\phi+C\int_{\mathbb{T}^{n}\times\mathbb{T}^{m}}\sigma^{p+1}dxd\phi\\
           &\leq Ck^2(p+1)^2\int_{\mathbb{T}^{n}\times\mathbb{T}^{m}} \sigma^{p+1}(1+\|D_xu\|^4)dxd\phi.
        \end{align*}
        Write $\alpha:=(1+\theta)^{1/2}>1$, then
        \begin{align*}
           \bigg(\int_{\mathbb{T}^{n}\times\mathbb{T}^{m}} \sigma^{(p+1)\alpha^2}dxd\phi\bigg)^{1/\alpha^2}
           &\leq Ck^2(p+1)^2\bigg(\int_{\mathbb{T}^{n}\times\mathbb{T}^{m}}\sigma^{(p+1)\alpha}dxd\phi\bigg)^{1/\alpha}
           \bigg(\int_{\mathbb{T}^{n}\times\mathbb{T}^{m}}(1+\|D_xu\|^4)^{\frac{\alpha}{\alpha-1}}dxd\phi\bigg)^{(\alpha-1)/\alpha}\\
           &\leq Ck^2(p+1)^2\bigg(\int_{\mathbb{T}^{n}\times\mathbb{T}^{m}}\sigma^{(p+1)\alpha}dxd\phi\bigg)^{1/\alpha},
        \end{align*}
        that is,
        \begin{align*}
          \bigg(\int_{\mathbb{T}^{n}\times\mathbb{T}^{m}} \sigma^{q\alpha}dxd\phi\bigg)^{1/\alpha}\leq Cq^{2\alpha}k^{2\alpha}\int_{\mathbb{T}^{n}\times\mathbb{T}^{m}}\sigma^qdxd\phi
        \end{align*}
        for each $q=(p+1)\alpha\geq\alpha$. Let $q=q_m:=\alpha^m$,
        \begin{align*}
           \bigg(\int_{\mathbb{T}^{n}\times\mathbb{T}^{m}} \sigma^{q_{m+1}}dxd\phi\bigg)^{1/q_{m+1}}
           \leq C^{1/q_m}q_m^{2\alpha/q_m}k^{2\alpha/q_m} \bigg(\int_{\mathbb{T}^{n}\times\mathbb{T}^{m}} \sigma^{q_m}dxd\phi\bigg)^{1/q_m},
        \end{align*}
        and iterating this process, we get
        \begin{align*}
           \|\sigma\|_{L^\infty(\mathbb{T}^{n}\times\mathbb{T}^{m})}
           \leq C^{\sum1/\alpha^m}\alpha^{\sum2m/(\alpha^m-1)}k^{\sum2/\alpha^{(m-1)}}
           \bigg(\int_{\mathbb{T}^{n}\times\mathbb{T}^{m}} \sigma^{\alpha}dxd\phi\bigg)^{1/\alpha}
           \leq Ck^{\frac{3\alpha-1}{\alpha-1}},
        \end{align*}
        where we used (\ref{2-A}) to estimate the last term of the above inequality.
        We restore the superscript $k$, because of $\sigma^k=e^{k[H(x,D_xu^k,\phi)-\bar{H}^k(P)]}$,
        \begin{align*}
           k[H(x,D_xu^k,\phi)-\bar{H}^k(P)]=\log \sigma^k\leq\log C+\frac{3\alpha-1}{\alpha-1}\log k,
        \end{align*}
        and then
        \begin{align}\label{4-3}
          \max_{(x,\phi)\in\mathbb{T}^{n}\times\mathbb{T}^{m}}H(x,D_xu^k,\phi)\leq\bar{H}^k(P)+\frac{C\log k}{k}.
        \end{align}
        Next, we are going to prove a kind of monotonicity. That is, for each $P\in\mathbb{R}^n$,  we have
        \begin{align}\label{4-4}
            \bar{H}^k(P)\leq \bar{H}^l(P),
        \end{align}
        if $k\leq l$. Indeed,
        \begin{align*}
           I_k[u^k]&=\int_{\mathbb{T}^{n}\times\mathbb{T}^{m}}e^{kH(x,D_xu^k,\phi)}dxd\phi
           \leq\int_{\mathbb{T}^{n}\times\mathbb{T}^{m}}e^{kH(x,D_xu^l,\phi)}dxd\phi\\
           &\leq\bigg(\int_{\mathbb{T}^{n}\times\mathbb{T}^{m}}e^{lH(x,D_xu^l,\phi)}dxd\phi\bigg)^{k/l}
           =(I_l[u^l])^{k/l},
        \end{align*}
        and take the log of both sides to derive (\ref{4-4}). Therefore, on account of $H$ grows quadratically in variable $p$, $(\ref{4-3})$ and $(\ref{4-4})$,
        \begin{align*}
          \max_{(x,\phi)\in\mathbb{T}^{n}\times\mathbb{T}^{m}} \bigg(\frac{\gamma}{2}\|D_xu^k\|^2-C\bigg)
          \leq \max_{(x,\phi)\in\mathbb{T}^{n}\times\mathbb{T}^{m}} H(x,D_xu^k,\phi)
          \leq\bar{H}^k(P)+\frac{C\log k}{k}
           \leq\bar{H}(P)+\frac{C\log k}{k},
        \end{align*}
        and then the conclusion holds.
     \end{proof}
      \subsection{\textbf{Proof of Lemma \ref{lem-4.1}}}
       \begin{proof}
             Let us first prove assertion (i).
             Recall that if $f(x)>0$, $x\in[a,b]$, then
             \begin{align*}
              \int_a^b f(x) dx \leq \log\bigg(\int_a^b e^{f(x)} dx\bigg).
              \end{align*}
              As an application, we have
               \begin{align*}
                \int_{\mathbb{T}^{n}\times\mathbb{T}^{m}}H(x,D_xu,\phi)dxd\phi
                 &\leq \displaystyle\mathop{\liminf}_{k\to\infty} \int_{\mathbb{T}^{n}\times\mathbb{T}^{m}} H(x,D_xu^k,\phi) dxd\phi\\
                 &\leq \displaystyle\mathop{\limsup}_{k\to\infty} \frac{1}{k}\log\bigg(\int_{\mathbb{T}^{n}\times\mathbb{T}^{m}}e^{kH(x,D_xu^k,\phi)}dxd\phi\bigg)\\
                 &\leq \bar{H}(P).
               \end{align*}
               So, $H(x,D_xu,\phi)\leq \bar{H}(P)$ almost everywhere in $\mathbb{T}^{n}\times\mathbb{T}^{m}$. Further,
               because Hamiltonian $H(x,y,\phi)$ is uniformly convex
              \begin{align}\label{A-1}
                  (D_xu-D_xu^\varepsilon)D_{yy}H\cdot(D_xu-D_xu^\varepsilon)\geq\frac{\gamma}{2}\|D_xu-D_xu^\varepsilon\|^2,
              \end{align}
              and then, by Taylor formula,
             \begin{align*}
                H(x,D_xu,\phi) &\geq H(x,D_xu^\varepsilon,\phi)
                +D_yH(x,D_xu^\varepsilon,\phi)\cdot(D_xu-D_xu^\varepsilon)+\frac{\gamma}{2}\|D_xu-D_xu^\varepsilon\|^2.
             \end{align*}
              Now, let us multiply both sides by $\rho_\varepsilon$  and integrate with respect to $\tilde{x}$  and $\tilde{\phi}$,
             \begin{align*}
                H(x,D_xu^\varepsilon,\phi)\leq\int_{\mathbb{R}^{n}\times\mathbb{R}^{m}} \rho_\varepsilon\bigg(H(x,D_xu,\phi)-\frac{\gamma}{2} |D_xu-D_xu^\varepsilon|^2\bigg)d\tilde{x}d\tilde{\phi}.
             \end{align*}
             Finally, based on $H(x,D_xu,\phi)\leq \bar{H}(P)$ almost everywhere in $\mathbb{T}^{n}\times\mathbb{T}^{m}$, we can get
            \begin{align}\label{A-8}
               H(x,D_xu^\varepsilon,\phi) \leq \bar{H}(P)+O(\varepsilon)-\chi_\varepsilon(x,\phi).
            \end{align}

             Next, we are start proving (ii).
             According to  \cite{L82,L83,L-P-V}, we can prove that there exists a unique constant $\bar{H}(P)$ such that
             $H(x,P+D_xv(x,P,\phi),\phi)=\bar{H}(P)$
             has a viscosity solution $v$.
             Then, based on  $v^k$ is minimizer of $I_k[\cdot]$, we can get
            \begin{align*}
              \int_{\mathbb{T}^{n}\times\mathbb{T}^{m}}e^{kH(x,P+D_xv^k,\phi)} dxd\phi
              \leq \int_{\mathbb{T}^{n}\times\mathbb{T}^{m}}e^{kH(x,P+D_xv,\phi)} dxd\phi
              =e^{k\bar{H}(P)},
            \end{align*}
           and consequently
           \begin{align*}
             \displaystyle\mathop{\limsup}_{k\to\infty}\frac{1}{k}\log\bigg(\int_{\mathbb{T}^{n}\times\mathbb{T}^{m}}
                                      e^{kH(x,P+D_xv^k,\phi)}dxd\phi\bigg)
             \leq\bar{H}(P).
           \end{align*}
             Assume next that, we take proof by contradiction, for some $\varepsilon>0$,
          \begin{align*}
             \displaystyle\mathop{\liminf}_{k\to\infty}\frac{1}{k}\log\bigg(\int_{\mathbb{T}^{n}\times\mathbb{T}^{m}}
                                      e^{kH(x,P+D_xv^k,\phi)}dxd\phi\bigg)
             <\bar{H}(P)-\varepsilon.
           \end{align*}
              Let $\Omega_{\varepsilon,k}:=\{(x,\phi)\in\mathbb{T}^{n}\times\mathbb{T}^{m}~|~H(x,D_xu^k,\phi) >\bar{H}(P)-\frac{\varepsilon}{2}\}$, for some sequence $\{k_j\}\to\infty$. We have
             \begin{align*}
               \bar{H}(P)-\varepsilon
               &\geq\frac{1}{k_j}\log\bigg(\int_{\mathbb{T}^{n}\times\mathbb{T}^{m}} e^{k_jH(x,P+D_xv^{k_j},\phi)}dxd\phi\bigg)
               \geq\frac{1}{k_j}\log\bigg(\int_{\Omega_{\varepsilon,k_j}} e^{k_jH(x,P+D_xv^{k_j},\phi)}dxd\phi\bigg)\\
                &\geq\frac{1}{k_j}\log\bigg(\int_{\Omega_{\varepsilon,k_j}} e^{k_j(\bar{H}(P)-\frac{\varepsilon}{2})}dxd\phi\bigg)
                =\frac{1}{k_j}\log|\Omega_{\varepsilon,k_j}|+\bar{H}(P)-\frac{\varepsilon}{2}.
              \end{align*}
             Hence
             \begin{align*}
              |\Omega_{\varepsilon,k_j}|\leq e^{-k_j\frac{\varepsilon}{2}},
              \end{align*}
              and by Fatou's lemma,
             \begin{align*}
               \int_{\mathbb{T}^{n}\times\mathbb{T}^{m}}H(x,D_xu,\phi)dxd\phi
               &\leq \displaystyle\mathop{\liminf}_{k_j\to\infty} \int_{\mathbb{T}^{n}\times\mathbb{T}^{m}} H(x,D_xu^{k_j},\phi)dxd\phi\\
                &=\displaystyle\mathop{\lim}_{k_j\to\infty} \int_{\mathbb{T}^{n}\times\mathbb{T}^{m}\backslash \Omega_{\varepsilon,k_j}}
                  +\int_{\mathbb{T}^{n}\times\mathbb{T}^{m}\cap \Omega_{\varepsilon,k_j}} H(x,D_xu^{k_j},\phi)dxd\phi
                \leq  \bar{H}(P)-\frac{\varepsilon}{2}.
              \end{align*}
             Further,
              \begin{align}\label{2-5}
              H(x,D_xu,\phi)\leq  \bar{H}(P)-\frac{\varepsilon}{2}
             \end{align}
              in $\mathbb{T}^{n}\times\mathbb{T}^{m}$.
             We set $u^\varepsilon=\rho_\varepsilon\ast u$ and due to (\ref{A-8}),
              \begin{align*}
              H(x,D_xu^\varepsilon,\phi)\leq  \bar{H}(P)-\frac{\varepsilon}{4}
             \end{align*}
             in $\mathbb{T}^{n}\times\mathbb{T}^{m}$,
             for sufficiently small $\varepsilon>0$.
             This contradiction can  be derived in terms of
             \begin{align*}
               \bar{H}(P)=\mathop {\inf}_{\hat{v}\in C^1(\mathbb{T}^{n})\times C(\mathbb{T}^{m})}\mathop {\max}_{(x,\phi)\in\mathbb{T}^{n}\times\mathbb{T}^{m}}H(x,P+D_x\hat{v},\phi).
             \end{align*}
            \end{proof}
\subsection{\textbf{Proof of Lemma \ref{lem-4.5}}}
\begin{proof}
        Because of the uniform convexity of the Hamiltonian  $H(x,y,\phi)$, the corresponding Lagrangian $L(x,\beta,\phi)$ also is uniformly convex with respect to $\beta$, then
       \begin{align}\label{3-3}
           \frac{\gamma}{2}\|\beta-r\|^2\leq L(x,\beta,\phi)-L(x,r,\phi)-D_\beta L(x,r,\phi)\cdot(\beta-r)
       \end{align}
       for each $\beta$, $r\in\mathbb{R}^n$, $x\in \mathbb{T}^n$, $\phi\in\mathbb{T}^m$.
       Let $r=D_yH(x,D_xu^\varepsilon,\phi)$ and integrate with respect $\mu$,
      \begin{align}\label{3-4}
         &\frac{\gamma}{2}\int_{\mathbb{R}^{m}\times\mathbb{T}^{n}\times\mathbb{T}^{m}}\|\beta-D_yH(x,D_xu^\varepsilon,\phi)\|^2d\mu\nonumber\\
         &\leq \int_{\mathbb{R}^{m}\times\mathbb{T}^{n}\times\mathbb{T}^{m}} L(x,\beta,\phi)-L(x,D_yH(x,D_xu^\varepsilon,\phi),\phi)-D_\beta L(x,D_yH(x,D_xu^\varepsilon,\phi),\phi)\cdot(\beta-D_yH(x,D_xu^\varepsilon,\phi))d\mu\nonumber\\
         &=\int_{\mathbb{R}^{m}\times\mathbb{T}^{n}\times\mathbb{T}^{m}} L(x,\beta,\phi) + H(x,D_xu^\varepsilon,\phi)-D_xu^\varepsilon\cdot D_yH(x,D_xu^\varepsilon,\phi)-D_xu^\varepsilon\cdot (\beta-D_yH(x,D_xu^\varepsilon,\phi)) d\mu\nonumber\\
         &=\int_{\mathbb{R}^{m}\times\mathbb{T}^{n}\times\mathbb{T}^{m}} L(x,\beta,\phi) + H(x,D_xu^\varepsilon,\phi)-D_xu^\varepsilon \cdot \beta d\mu.
      \end{align}
      According to Lemma \ref{lem-4.1} $(ii)$,
      \begin{align}\label{3-5}
         \int_{\mathbb{T}^{n}\times\mathbb{T}^{m}} \chi_\varepsilon(x,\phi)d\sigma\leq \bar{H}(P)+O(\varepsilon)
         -\int_{\mathbb{T}^{n}\times\mathbb{T}^{m}} H(x,D_xu^\varepsilon,\phi) d\sigma.
      \end{align}
      Combining (\ref{3-4}) and (\ref{3-5}), we have
      \begin{align}\label{3-6}
          &\frac{\gamma}{2}\int_{\mathbb{R}^{m}\times\mathbb{T}^{n}\times\mathbb{T}^{m}}\|\beta-D_yH(x,D_xu^\varepsilon,\phi)\|^2d\mu
          +\int_{\mathbb{T}^{n}\times\mathbb{T}^{m}} \chi_\varepsilon(x,\phi)d\sigma\nonumber\\
          &\leq \int_{\mathbb{R}^{m}\times\mathbb{T}^{n}\times\mathbb{T}^{m}}L(x,\beta,\phi)-D_xu^\varepsilon\cdot\beta d\mu+\bar{H}(P)+O(\varepsilon)\nonumber\\
          &=\bar{L}(Q)-P\cdot Q+\bar{H}(P)+O(\varepsilon)=O(\varepsilon).
      \end{align}
      Thereby, by H$\ddot{o}$der inequality,
      \begin{align*}
          \int_{\mathbb{T}^{n}\times\mathbb{T}^{m}} \chi^{1/2}_\varepsilon(x,\phi)d\sigma
          \leq \bigg(\int_{\mathbb{T}^{n}\times\mathbb{T}^{m}} \chi_\varepsilon(x,\phi)d\sigma\bigg)^{1/2}\leq C\varepsilon^{1/2},
      \end{align*}
      let $\varepsilon=1/2^k$,
      \begin{align*}
          \int_{\mathbb{T}^{n}\times\mathbb{T}^{m}} \sum^{\infty}_{k=1} \chi^{1/2}_{1/2^k}(x,\phi) d\sigma\leq C \sum^{\infty}_{k=1} (1/2^k)^{1/2}<\infty.
      \end{align*}
      So,
      \begin{align*}
          \sum^{\infty}_{k=1} \chi^{1/2}_{1/2^k}(x,\phi)<\infty
      \end{align*}
      for $\sigma$-almost everywhere point $(x,\phi)\in\mathbb{T}^{n}\times\mathbb{T}^{m}$.
      Fix such a point $(x,\phi)$, for some constant $\kappa>0$, we have
      \begin{align*}
          \chi^{1/2}_{2r}(x,\phi)
          &\geq\kappa\bigg(\frac{1}{|B(x,r)\times B(\phi,r)|}\int_{B(x,r)\times B(\phi,r)}\|D_xu(\tilde{x},\tilde{\phi})-(D_xu)_r\|^2d\tilde{x}d\tilde{\phi}\bigg)^{1/2}\\
          &\geq\kappa\bigg(\frac{1}{|B(x,r)\times B(\phi,r)|}\int_{B(x,r)\times B(\phi,r)}\|D_xu(\tilde{x},\tilde{\phi})-(D_xu)_r\|d\tilde{x}d\tilde{\phi}\bigg),
      \end{align*}
      where
      \begin{align*}
          (D_xu)_r=\frac{1}{|B(x,r)\times B(\phi,r)|}\int_{B(x,r)\times B(\phi,r)}D_xu(\tilde{x},\tilde{\phi}) ~d\tilde{x}d\tilde{\phi}.
      \end{align*}
      Hence
      \begin{align*}
          \sum^\infty_{k=1} \frac{1}{|B(x,1/2^{k+1})\times B(\phi,1/2^{k+1})|}
          &\int_{B(x,1/2^{k+1})\times B(\phi,1/2^{k+1})}\|D_xu(\tilde{x},\tilde{\phi})-(D_xu)_{1/2^{k+1}}\|~d\tilde{x}d\tilde{\phi} \\
          &\leq \frac{1}{\kappa} \sum^{\infty}_{k=1}\chi^{1/2}_{1/2^k} < \infty.
      \end{align*}
      Since
       \begin{align*}
           &|(D_xu)_{1/2^{k+1}}-(D_xu)_{1/2^{k}}|\\
           &= \bigg|\frac{1}{|B(x,1/2^{k+1})\times B(\phi,1/2^{k+1})|} \int_{B(x,1/2^{k+1})\times B(\phi,1/2^{k+1})} D_xu(\tilde{x},\tilde{\phi})-(D_xu)_{1/2^{k}} d\tilde{x}d\tilde{\phi} \bigg|\\
           &\leq \frac{1}{|B(x,1/2^{k+1})\times B(\phi,1/2^{k+1})|} \int_{B(x,1/2^{k+1})\times B(\phi,1/2^{k+1})} |D_xu(\tilde{x},\tilde{\phi})-(D_xu)_{1/2^{k}}| d\tilde{x}d\tilde{\phi}\\
           &\leq C\frac{1}{|B(x,1/2^{k})\times B(\phi,1/2^{k})|} \int_{B(x,1/2^{k})\times B(\phi,1/2^{k})} |D_xu(\tilde{x},\tilde{\phi})-(D_xu)_{1/2^{k}}| d\tilde{x}d\tilde{\phi},
       \end{align*}
       we obtain
       \begin{align*}
          \sum^{\infty}_{k=1} |(D_xu)_{1/2^{k+1}}-(D_xu)_{1/2^{k}}| <\infty,
       \end{align*}
        and hence the limit $D_xu(x,\phi):=\displaystyle \mathop{\lim}_{k \to \infty}(D_xu)_{(x,\phi),1/2^{k}}$ exists.
        So, $u$ is differentiable at $x$, and $\chi_\varepsilon\rightarrow0$, as $\varepsilon\rightarrow0$.
        Therefore,
        \begin{align*}
            D_xu^\varepsilon\rightarrow D_xu,
        \end{align*}
       $\sigma$-almost everywhere, and
       \begin{align*}
            \int_{\mathbb{R}^{m}\times\mathbb{T}^{n}\times\mathbb{T}^{m}}\|\beta-D_yH(x,D_xu,\phi)\|^2d\mu=0.
       \end{align*}
       In particularly,
       \begin{align*}
          Q=\int_{\mathbb{T}^{n}\times\mathbb{T}^{m}} D_yH(x,D_xu,\phi)d\sigma.
       \end{align*}

        \end{proof}
       \section*{Acknowledgment}
			The work of K. Wang is partially  supported by NSFC (Nos. 12171315, 11931016).
The work of Y. Li is partially supported by National Basic Research Program of China (No. 2013CB834100), NSFC (No. 12071175), Special Funds of Provincial Industrial Innovation of Jilin Province China (Grant No. 2017C028-1), Project of Science and Technology Development of Jilin Province China (Grant No. 20190201302JC).
	

\section{References}
		\begin{thebibliography}{99}
           \bibitem{A-S-V} A. Arapostathis, S. S Sastry, P. P. Varaiya,  Global analysis of swing dynamics. IEEE Trans. Circuits and Systems 29 (1982), no. 10, 673-679.
        \bibitem{C-I-P-P} G. Contreras, R. Iturriaga, G. P. Paternain, M. Paternain, Lagrangian graphs, minimizing measures and Ma\~{n}\'{e}'s critical values. Geom. Funct. Anal. 8 (1998), no. 5, 788-809.
           \bibitem{E98}  L. C. Evans, Partial differential equations. American Mathematical Society, Providence, RI, 1998.
           \bibitem{E02} L. C. Evans, D. Gomes, Effective Hamiltonians and averaging for Hamiltonian dynamics. II. Arch. Ration. Mech. Anal. 161 (2002), no.4, 271-305.
        \bibitem{E03} L. C. Evans,   Some new PDE methods for weak KAM theory. Calc. Var. Partial Differential Equations 17 (2003), no. 2, 159-177.
        \bibitem{E04} L. C. Evans, Towards a quantum analog of weak KAM theory. Comm. Math. Phys. 244 (2004), no. 2, 311-334.

        \bibitem{F97} A. Fathi, Th\'{e}or\`{e}me KAM faible et th\'{e}orie de Mather sur les syst\`{e}mes lagrangiens. C. R. Acad. Sci. Paris
           S\'{e}r. I Math. 324 (1997), no. 9, 1043-1046.
        \bibitem{F98} A. Fathi, Sur la convergence du semi-groupe de Lax-Oleinik. C. R. Acad. Sci. Paris S\'{e}r. I Math. 327 (1998), no.3,  267-270.
           \bibitem{G-S} A. Gholami, X. Sun, The impact of damping in second-order dynamical systems with applications to power grid stability.  SIAM J. Appl. Dyn. Syst. 21 (2022), no. 1, 405-437.
             \bibitem{G-T}  D. Gilbarg, N. S. Trudinger,  Elliptic partial differential equations of second order. Springer-Verlag, Berlin, 1983.
        \bibitem{G} D. Gomes,  Perturbation theory for viscosity solutions of Hamilton-Jacobi equations and stability of Aubry-Mather sets. SIAM J. Math. Anal. 35 (2003), no. 1, 135-147.
            \bibitem{K} H. K. Khalil,  Nonlinear systems. Macmillan Publishing Company, New York, 1992.
       \bibitem{L}   G. M. Lieberman,  Solvability of quasilinear elliptic equations with nonlinear boundary conditions. Trans. Amer. Math. Soc. 273 (1982), no. 2, 753-765.
        \bibitem{L82} P. L.  Lions,  Generalized solutions of Hamilton-Jacobi equations. Research Notes in Mathematics, 69. Pitman, London, 1982.
        \bibitem{L83} P. L. Lions,   Existence results for first-order Hamilton-Jacobi equations. Ricerche Mat. 32 (1983), no. 1, 3-23.
        \bibitem{L-P-V} P. L. Lions, G. Papanicolaou, S. Varadhan,   Homogenization of Hamilton-Jacobi equation. unpublished preprint, 1987.
          \bibitem{M} J. N. Mather, Action minimizing invariant measures for positive definite Lagrangian systems,
Math. Z. 207 (1991), 169-207.
        \bibitem{N-J-L} X. Niu, S. Ji, Y. Li, Relative Equilibrium Via Viscosity Solution.  J. Differential Equations 392 (2024), 325-363.
        \bibitem{N-W-L} X. Niu, K. Wang, Y. Li, Resonance conjecture via weak KAM theory. J. Math. Pures Appl. (9) 172 (2023), 139-163.
          \bibitem{Q-M-K-Z}   Q. Qiu, R. Ma, J. Kurths, M. Zhan,  Swing equation in power systems: approximate analytical solution and bifurcation curve estimate. Chaos 30 (2020), no. 1,  11 pp.

          \bibitem{S-M-V}   F. M. A.  Salam, J. E. Marsden, P. P. Varaiya,  Arnol'd diffusion in the swing equations of a power system. IEEE Trans. Circuits and Systems 31 (1984), no. 8, 673-688.
              \bibitem{S-G-H}  T. H.  Scholl, L. Gr\"{o}ll, V. Hagenmeyer,  Time delay in the swing equation: a variety of bifurcations. Chaos 29 (2019), no. 12,  12 pp.
            \bibitem{S} S. J. Skar, Stability of multimachine power systems with nontrivial transfer conductances.
SIAM J. Appl. Math. 39 (1980), no. 3, 475-491.
          \bibitem{S-M-H}  Y. Susuki, I. Mezi\'{c}, T. Hikihara,  Coherent swing instability of power grids. J. Nonlinear Sci. 21 (2011), no. 3, 403-439.
        \bibitem{W-L}  K. Wang, Y. Li,  Some results on weak KAM theory for time-periodic Tonelli Lagrangian systems. Adv. Nonlinear Stud. 13 (2013), no. 4, 853-866.
        \bibitem{W-Y}K.  Wang, J. Yan,  A new kind of Lax-Oleinik type operator with parameters for time-periodic positive definite Lagrangian systems. Comm. Math. Phys. 309 (2012), no. 3, 663-691.
         \bibitem{X}C. Xu, Regularity for quasilinear second-order subelliptic equations. Comm. Pure Appl. Math. 45 (1992), no. 1, 77-96.
        \end{thebibliography}
\end{document}